\documentclass[a4paper,oneside]{article}
\pdfoutput=1

\RequirePackage{amsmath, amsthm, amssymb}
\RequirePackage{mathtools}
\RequirePackage{natbib}
\RequirePackage{tikz}

\RequirePackage[margin=1in]{geometry}
\RequirePackage{libertine}
\RequirePackage{courier}
\RequirePackage{mathpazo}
\RequirePackage{xcolor}
\RequirePackage{pifont}
\RequirePackage{authblk}
\RequirePackage[backref=page]{hyperref}
\RequirePackage{enumerate}
\RequirePackage{titlesec}

\definecolor{wred}{rgb}{0.533333,0.10980,0.15686}
\definecolor{wredlight}{rgb}{0.788, 0.11, 0.157}
\titleformat{\section}
{\color{wred}\normalfont\Large\bfseries}
{\color{wred}\thesection}{1em}{}
\titleformat{\subsection}
{\color{wred}\normalfont\Large\bfseries}
{\color{wred}\thesubsection}{1em}{}

\hypersetup{
	colorlinks=true,
  linkcolor=wredlight,
  citecolor=wredlight,
  filecolor=wredlight,
  urlcolor=wredlight
}

 \renewcommand*{\backrefalt}[4]{%
    \ifcase #1%
     \or (page:~#2)%
     \else (pages:~#2)%
    \fi%
    }

    \makeatletter
\def\@fnsymbol#1{\ensuremath{\ifcase#1\or \dagger\or \ddagger\or
   \mathsection\or \mathparagraph\or \|\or **\or \dagger\dagger
   \or \ddagger\ddagger \else\@ctrerr\fi}}
    \makeatother

\newtheorem{theorem}{Theorem}[section]
\newtheorem{lemma}{Lemma}[section]
\newtheorem{definition}{Definition}[section]

\newtheorem{example}{Example}[section]
\newtheorem{corollary}{Corollary}[section]
\newtheorem{remark}{Remark}[section]

\providecommand{\keywords}[1]
{
  \small	
  \textbf{\textit{Keywords---}} #1
}
\definecolor{darkgreen}{rgb}{0.2, 0.5, 0.2}

\DeclareMathOperator{\diag}{diag}
\DeclareMathOperator{\Trace}{Tr}
\DeclareMathOperator*{\argmin}{arg\,min}

\newcommand{\eqdef}{\doteq}
\newcommand{\pred}[1]{\delta\left[#1\right]}
\newcommand{\N}{\mathbb{N}}
\newcommand{\trn}{^\intercal}
\newcommand{\abs}[1]{\left| #1 \right|}
\newcommand{\nrm}[1]{\left\Vert #1 \right\Vert}
\newcommand{\PR}[2][]{\mathbb{P}_{#1}\left( #2 \right)}
\newcommand{\E}[2][]{\mathbb{E}_{#1}\left[ #2 \right]}
\newcommand{\Var}[2][]{\mathbb{V}_{#1}\left[ #2 \right]}
\newcommand{\eps}{\varepsilon}
\newcommand{\bigO}{\mathcal{O}}
\newcommand{\set}[1]{\left\{ #1 \right\}}
\newcommand{\R}{\mathbb{R}}

\newcommand{\calX}{\mathcal{X}}

\newcommand{\calP}{\mathcal{P}}

\newcommand{\calF}{\mathcal{F}}
\newcommand{\calB}{\mathcal{B}}
\newcommand{\calM}{\mathcal{M}}

\newcommand{\calH}{\mathcal{H}}

\newcommand{\calN}{\mathcal{N}}

\newcommand{\calK}{\mathcal{K}}
\newcommand{\calA}{\mathcal{A}}
\newcommand{\bbP}{\mathbb{P}}
\newcommand{\bbN}{\mathbb{N}}

\newcommand{\bbR}{\mathbb{R}}
\newcommand{\bbQ}{\mathbb{Q}}
\newcommand{\bbH}{\mathbb{H}}
\newcommand{\bbE}{\mathbb{E}}
\newcommand{\bbV}{\mathbb{V}}

\newcommand{\tmix}{t_{\mathsf{mix}}}
\newcommand{\betamix}{t^{\beta}_{\mathsf{mix}}}
\newcommand{\kme}{\mu_{\mathbb{P}}}
\newcommand{\ekme}{\widehat{\mu}_{\mathbb{P}}}
\newcommand{\tv}[1]{\nrm{#1}_{\mathsf{TV}}}

\author[1]{Geoffrey Wolfer \thanks{email: geo.wolfer@aoni.waseda.jp. \\
Part of this research was supported by the Special Postdoctoral Researcher (SPDR) program of RIKEN.}}
\author[2]{Pierre Alquier \thanks{email: alquier@essec.edu.}}
\affil[1]{Waseda University, Center for Data Science}
\affil[2]{ESSEC Business School, Asia-Pacific campus}

\date{\date{}}

\title{\vspace{-1.0cm} Variance-Aware Estimation of \\ Kernel Mean Embedding}

\begin{document}
\maketitle

\begin{abstract}
    An important feature of kernel mean embeddings (KME) is that the rate of convergence of the empirical KME to the true distribution KME can be bounded independently of the dimension of the space, properties of the distribution and smoothness features of the kernel. We show how to speed-up convergence by leveraging variance information in the reproducing kernel Hilbert space. Furthermore, we show that even when such information is a priori unknown, we can efficiently estimate it from the data, recovering the desiderata of a distribution agnostic bound that enjoys acceleration in fortuitous settings. We further extend our results from independent data to stationary mixing sequences and illustrate our methods in the context of hypothesis testing and robust parametric estimation.
\end{abstract}

\keywords{kernel mean embedding, maximum mean discrepancy, distribution learning, empirical Bernstein}

\tableofcontents

\section{Introduction}
\label{section:introduction}
Estimating a probability distribution $\bbP$ over a space $\calX$ from $n$ iid samples $X_1, \dots, X_n \sim \bbP$ 
is a central problem in computer science and statistics \citep{kearns1994learnability, tsybakov2008introduction}.
To formalize the question, one selects a distance (or at least a contrast function) between distributions and oftentimes introduces assumptions on the underlying probability space (for instance finitely supported, probability function is absolutely continuous with respect to the Lebesgue measure, H\"older continuous, \dots).
Increasing the stringency of assumptions generally leads to substantially faster minimax rates. 
For instance, in the finite support $\abs{\calX} < \infty$ case and with respect to total variation, whilst the expectation risk evolves roughly as $\sqrt{\abs{\calX} / n}$, it is known that this rate can be sharpened by replacing $\abs{\calX}$ with a bound on the ``half-norm" \footnote{The half-norm is not a norm in the proper sense as it violates the triangle inequality.} $\|\bbP\|_{1/2} \eqdef (\sum_{x \in \calX} \sqrt{\bbP(x)})^2$  \citep{berend2013sharp}, when defined, and which corresponds to some measure of entropy of the underlying distribution. Furthermore, even without prior knowledge of $\|\bbP\|_{1/2}$, one can construct confidence intervals around $\bbP$ that depend on $\|\widehat{\bbP}_n\|_{1/2}$---the half-norm of the empirical distribution $\widehat{\bbP}_n(x) = \frac{1}{n} \sum_{t=1}^{n} \pred{X_t = x}$ \citep[Theorem~2.1]{ckw2020}. Namely, when $\bbP$ is supported on $\N$, with probability at least $1 - \delta$, it holds that
\begin{equation}
\label{eq:tv-confidence-interval}
\begin{split}
    \tv{\widehat{\bbP}_n - \bbP} \leq \frac{1}{\sqrt{n}} \left(  \nrm{\widehat{\bbP}_n}_{1/2}^{1/2} + 3 \sqrt{\frac{1}{2} \log (2/\delta)} \right).
\end{split}
\end{equation}
The advantage of the above expression is twofold $(a)$ we do not need make assumptions on $\|\bbP\|_{1/2}$, which could be non-convergent\footnote{As opposed to the empirical proxy $\|\widehat{\bbP}_n\|_{1/2}$, which is always a finite quantity.}, and $(b)$ in favorable cases where $\|\widehat{\bbP}_n\|_{1/2}$ is small, the intervals will be narrower.
In this paper, we set out to explore the question of whether analogues of \eqref{eq:tv-confidence-interval} are possible for general probability spaces and with respect to maximum mean discrepancy.

\subsection{Notation and Background}
The set of integers up to $n \in \N$ is denoted by $[n] \eqdef \set{1, 2, \dots}$.  Let $\calX$ be a separable topological space, and $\calP(\calX)$ the set of all Borel probability measures over $\calX$. 
For a bounded function $f\colon \calX \to \R $, we define for convenience
\begin{equation}
\label{definition:shorthands}
    \overline{f} \eqdef \sup_{x \in \calX} f(x) \qquad \underline{f} \eqdef \inf_{x \in \calX} f(x), \qquad \Delta f \eqdef \overline{f} - \underline{f}.
\end{equation}
Let $k \colon \calX \times \calX \to \R$ be a continuous positive definite  kernel and $\calH_k$ be the associated reproducing kernel Hilbert space (RKHS) \citep{berlinet2011reproducing}.
We assume that the kernel is bounded\footnote{Note that the supremum of the kernel is always reached on the diagonal, that is, $\sup_{(x,x') \in \calX^2} k(x,x')=\sup_{x \in \calX} k(x,x)$.} in the sense that $\sup_{x \in \calX} k(x,x) < \infty$. Letting $\bbP \in \calP(\calX)$, the kernel mean embedding (KME) of $\bbP$ is defined as
\begin{equation*}
    \kme \eqdef \E[X \sim \bbP]{k(X, \cdot)} = \int_{\calX} k(x, \cdot) d\bbP(x) \in \calH_k,
\end{equation*}
which is interpreted as a Bochner integral \citep[Chapter~2]{diestel1977vector}.
Let $n \in \N$ and $X_1, \dots, X_n$ be a sequence of observations sampled independently\footnote{Unless otherwise specified, the data will be considered iid. We will address the case of dependent data in Section~\ref{section:dependent-data}.} from $\bbP$. We write
\begin{equation}
\label{eq:kme-empirical-distribution}
    \qquad
    \ekme(X_1, \dots, X_n) = 
    \ekme \eqdef \frac{1}{n} \sum_{t = 1}^{n} k(X_t, \cdot) \in \calH_k,
\end{equation}
for the KME of the empirical distribution $\widehat{\bbP}_n \eqdef \frac{1}{n} \sum_{t=1}^{n}\delta_{X_t}$,
and call the distance $\nrm{\ekme - \kme}_{\calH_k}$ the
maximum
mean discrepancy (MMD) between the true mean embedding and its empirical estimator.
A kernel is called characteristic if the map $\mu \colon \bbP \mapsto \mu_{\bbP}$ is injective. This property ensures that $\nrm{\mu_{\bbP} - \mu_{\bbP'}}_{\calH_k} = 0 $ if and only if $\bbP = \bbP'$.
A kernel is called translation invariant\footnote{Also sometimes referred to as an \emph{anisotropic stationary} kernel.} (TI) when there exists a positive definite function $\psi \colon \calX \to \R$ such that for any $x,x' \in \calX$, 
\begin{equation}
\label{definition:translation-invariant}
    k(x,x') = \psi(x' - x).
\end{equation}
In particular, $\overline{\psi} = \psi(0) = \overline{k}$.
When $\calX = \R^d$, a kernel $k$ is said to be a radial basis function (RBF) when
for any $x,x' \in \calX$, $k(x,x') = \phi(\nrm{x' - x}_2)$ for some function $\phi \colon \R_+ \to \R$. 
Noticeably, $k$ being positive definite does not preclude it from taking negative values.
However, when $k$ is RBF, the following lower bound on $\underline{\phi}$ holds (see for instance \citealt{genton2001classes, stein1999interpolation}),
\begin{equation*}
    \underline{\phi}  \geq \overline{\phi} \inf_{t \geq 0} \set{ \left( \frac{2}{t}\right)^{(d-2)/2} \Gamma(d/2) J_{(d-2)/2}(t) },
\end{equation*}
where $\Gamma$ is the Gamma function and $J_\beta$ is the Bessel function of the first kind of order $\beta$, 
showing that $|\underline{\phi}|$ 
becomes evanescent as the dimension increases.

\subsection{Related Work}
From an asymptotic standpoint, the weak law of large numbers asserts that $\ekme$ converges to the true $\kme$.
Furthermore,
$$\sqrt{n}\left(\ekme - \kme\right)$$
converges in distribution to a zero mean Gaussian process on $\calH_k$ \citep[Section~9.1]{berlinet2011reproducing}.
This work, however, is more concerned with the finite sample theory, and more specifically with the rate of convergence of $\ekme$ towards $\kme$ with respect to the RKHS norm.
Conveniently and perhaps surprisingly, it 
 is possible to derive a rate that depends
 neither on the smoothness of the considered kernel $k$, nor the properties of the true distribution.
To obtain a distribution independent rate at $\bigO_P(n^{-1/2})$, a typical strategy (see for example \citealt[Section~B.1]{lopez2015towards}) consists in first expressing the dual relationship between the norm in the RKHS and the uniform norm of an empirical process,
\begin{equation}
\label{eq:dual-relationship-rkhs-norm-empirical-process}
    \nrm{\ekme - \kme}_{\calH_k} = \sup_{\substack{f \in \calH_k \\ \nrm{f}_{\calH_k} \leq 1}} \langle f, \ekme - \kme  \rangle_{\calH_k} = 
    \sup_{\substack{f \in \calH_k \\ \nrm{f}_{\calH_k} \leq 1}} \left(\frac{1}{n} \sum_{t = 1}^{n} f(X_t) - \bbE f \right).
\end{equation}
A classical symmetrization argument \citep{mohri2018foundations} followed by an application of McDiarmid's inequality  \citep{mcdiarmid1989method} yield that with probability at least $1 - \delta$,
\begin{equation}
\label{eq:rademacher-complexity-control}
    \sup_{\substack{f \in \calH_k \\ \nrm{f}_{\calH_k} \leq 1}} \left(\frac{1}{n} \sum_{t = 1}^{n} f(X_t) - \bbE f \right) \leq 2 \mathfrak{R}_n + \sqrt{\frac{2 \overline{k} \log (1/\delta)}{n}},
\end{equation}
where $\mathfrak{R}_n$ is the Rademacher complexity \citep[Definition~3.1]{mohri2018foundations} of the class of unit functions in the RKHS.
The bound
$\mathfrak{R}_n^2
    \leq \overline{k}/n$ \citep{bartlett2002rademacher},
and an application of Jensen's inequality conclude the claim.
With a more careful analysis, \citet[Proposition~A.1]{tolstikhin2017minimax} halved the constant of the first term in \eqref{eq:rademacher-complexity-control},
hence showed that with probability at least $1 - \delta$,
\begin{equation}
\label{eq:tolstikhin-bound}
    \nrm{\ekme - \kme}_{\calH_k} \leq  \sqrt{\frac{\overline{k}}{n}} + \sqrt{\frac{2 \overline{k} \log (1/ \delta)}{n}}.
\end{equation}
What is more, \citet[Theorems~1,6,8]{tolstikhin2017minimax} provide corresponding lower bounds in $\Omega_P(n^{-1/2})$, 
showing that the embedding of the empirical measure achieves an unimprovable rate of convergence. Results similar to~\eqref{eq:tolstikhin-bound} (with worse constants) can be derived when the observations are not independent, see \citet[Lemma 7.2]{cherief2022finite}.
We note that other estimators have been proposed in the literature for the case where more information is available about the underlying distribution $\bbP$.
See for instance \citet{muandet2014kernel, muandet2016kernel} who propose a shrinkage estimator inspired by the James-Stein estimator.
In this work, we stress that our estimator will be the KME of the empirical distribution, defined in \eqref{eq:kme-empirical-distribution}.

\subsection{Main Contributions}
In this paragraph, we briefly summarize our findings.
The tilde notation suppresses constants and logarithmic factors in $1/\delta$, and the reader is referred to Theorem~\ref{theorem:variance-aware-mmd-estimation} and Theorem~\ref{theorem:empirical-variance-mmd-bound} for the full expression of the second order terms.
Our first contribution (Theorem~\ref{theorem:variance-aware-mmd-estimation}) consists in deriving
a bound on $\nrm{\ekme - \kme}_{\calH_k}$ that can leverage additional information about the underlying distribution $\bbP$ and the selected kernel $k$.
Namely, we show that with probability at least $1 - \delta$,
\begin{equation*}
    \nrm{\ekme - \kme}_{\calH_k}  \leq \sqrt{2 v_k(\bbP) \frac{\log (2/ \delta)}{n}} + \widetilde{\bigO} \left( \frac{1}{n} \right),
\end{equation*}
where $v_k(\bbP)$---defined in \eqref{eq:rkhs-variance}---corresponds to some notion of variance in the RKHS. Notably, it holds that
 $v_k(\bbP) \leq \overline{k}$,
hence our upper bound is superior to and able to recover known bounds in most settings (Remark~\ref{remark:superior-most-settings}).

Our chief technical contribution (Theorem~\ref{theorem:empirical-variance-mmd-bound}) is in establishing that, at least when $k$ is translation invariant, 
the dependence of the above confidence interval can be made independent of the underlying distribution by replacing the variance by an empirical proxy.
Specifically, with probability at least $1- \delta$,
\begin{equation*}
\begin{split}
    \nrm{\ekme - \kme}_{\calH_{k}} &\leq  \sqrt{2 \widehat{v}_k(X_1, \dots, X_n) \frac{ \log (4/ \delta)}{n}} + \widetilde{\bigO}\left( \frac{1}{n} \right) \\
    \end{split}
    \end{equation*}
where $\widehat{v}_k$---defined in \eqref{eq:general-empirical-variance-proxy},  \eqref{eq:translation-invariant-empirical-variance-proxy}---is a quantity that only depends on the chosen kernel and the sample. We also extend the latter result to kernels without the translation invariance property in
Theorem~\ref{theorem:general-empirical-variance-mmd-bound}.

Expanding beyond the iid setting, we obtain convergence rates for stationary mixing sequences. For $\phi$-mixing processes, Theorem~\ref{theorem:variance-aware-mmd-estimation-phi-mixing} establishes the following rate of convergence.
\begin{equation*}
    \nrm{\ekme - \kme}_{\calH_{k}}
     \leq  \sqrt{\frac{v_k + \Sigma_n}{n}} + 4 \sqrt{\frac{2 v_k \nrm{\Gamma}_2 \log (1 / \delta)}{n}} + \widetilde{\bigO}\left( \frac{1}{n} \right),
\end{equation*}
where $\Sigma_n$ is a total measure of covariance between observations in the RKHS (refer to Lemma~\ref{lemma:exact-expression-second-moment-covariances}) and $\nrm{\Gamma}_2$ is the spectral norm of a coupling matrix defined in \eqref{definition:coupling-matrix}. We additionally analyze the broader class of $\beta$-mixing processes in Theorem~\ref{theorem:variance-aware-mmd-estimation-beta-mixing}.

\subsection{Outline}
In Section~\ref{section:variance-aware}, for the problem of estimating a distribution with respect to maximum mean discrepancy, we give a convergence rate with a dominant term in $\bigO_P(n^{-1/2})$ involving a variance term $v$ which depends on both the chosen kernel $k$ and the underlying distribution $\bbP$.
We then illustrate how the rate can subdue minimax lower bounds when $v$ is favorable.
In particular, we show that for large collections of kernels, and when $\calX \subset \R^d$, the variance  
$v$ can be controlled by a quantity that decouples the influence of the kernel and a measure of total variance of $\bbP$.
In Section~\ref{section:variance-empirical}, we proceed to show that even if $v$ is unknown, it is possible to efficiently estimate it from the data with an ``empirical Bernstein" approach whenever the kernel is translation invariant, or at least enjoys a mildly varying diagonal.
In Section~\ref{section:two-means}, we illustrate the benefits of the variance aware bounds on the classical two-sample test problem.
In Section~\ref{section:dependent-data}, we establish convergence rates for stationary $\phi$ and $\beta$ mixing sequences.
In Section~\ref{section:parametric-estimation}, we put our methods into practice, first in the context of hypothesis testing, and second by improving the results of \cite{briol2019statistical} and \citet{cherief2022finite} in the context of robust parametric maximum mean discrepancy estimation.
All proofs are deferred to Section~\ref{section:proofs} for clarity of the exposition.

\section{Variance-Aware Convergence Rates}
\label{section:variance-aware}
The central quantity we propose to consider is the following variance term in the RKHS,
\begin{equation}
\label{eq:rkhs-variance}
    v_k(\bbP) \eqdef \bbE_{X \sim \bbP} \nrm{k(X, \cdot) - \kme}_{\calH_k}^2.
\end{equation}
It is clear that $v_k(\bbP)$ depends both on the choice of kernel $k$ and on the underlying distribution, and we will simply write $v = v_k(\bbP)$ to avoid encumbering notation.
Simple calculations (refer to Lemma~\ref{lemma:exact-expression-second-moment-covariances}) show that
$$\bbE \nrm{\ekme - \kme}_{\calH_k}^2 = \frac{v}{n},$$
where the expectation is taken with respect to an independent sample $X_1, \dots, X_n \sim \bbP$,
thus by applications of Jensen's and Chebyshev's inequalities, we readily obtain a rate on the expected risk in terms of $v$,
\begin{equation*}
    \begin{split}
        \bbE  \nrm{\ekme - \kme}_{\calH_k} &\leq \sqrt{\frac{v}{n}},
    \end{split}
\end{equation*}
and a deviation bound
\begin{equation*}
    \begin{split}
        \PR{ \nrm{\ekme - \kme}_{\calH_k} > (1 + \tau) \sqrt{v/n}} &\leq 1/(1 + \tau)^2, \tau \in (0, \infty).
    \end{split}
\end{equation*}
One basic feature of RKHS is that new kernels can be constructed by linearly combining existing kernels---see for instance \citet{gretton2013introduction}. It is therefore noteworthy that $v$ is linear with respect reproducing kernels.
Reformulate
\begin{equation*}
    v = \E[X \sim \bbP]{k(X,X)} - \E[X,X' \sim \bbP]{k(X,X')},
\end{equation*}
and suppose that $k = \sum_{k_i \in \calK} \alpha_i k_i$ where $\calK = \set{k_1, \dots, k_{\abs{\calK}}}$ is a collection of reproducing kernels with $\alpha \in \R^{\abs{\calK}}$. It then holds that
\begin{equation*}
    v_{k} = \sum_{k_i \in \calK} \alpha_i v_{ k_i}.
\end{equation*}
Moving on to high-probability confidence bounds,
an application of Bernstein's inequality in Hilbert spaces  \citep{pinelis1985remarks,yurinsky1995sums} not only recovers the rate of convergence, as alluded to by \citet{tolstikhin2017minimax}, but in fact yields a maximal inequality.

\begin{theorem}[Variance-aware confidence interval]
\label{theorem:variance-aware-mmd-estimation}
Let $X_1, \dots, X_n \stackrel{iid}{\sim} \bbP$, $k \colon \calX \times \calX \to \R$ be a reproducing kernel, and  let
\begin{equation*}
    v \eqdef \bbE_{X \sim \bbP} \nrm{k(X, \cdot) - \kme}_{\calH_k}^2.
\end{equation*}
With probability at least $1 - \delta$, it holds that
\begin{equation*}
    \max_{1 \leq t \leq n} \set{ t
    \nrm{\ekme(X_1, \dots, X_t) - \kme}_{\calH_k} } \leq  \sqrt{2 v n \log (2/ \delta)} + \frac{4}{3}\sqrt{\overline{k}}\log (2/ \delta),
\end{equation*}
where $\overline{k} = \sup_{x \in\calX} k(x,x) $.
In particular,
with probability at least $1 - \delta$, it holds that
\begin{equation*}
    \nrm{\ekme(X_1, \dots, X_n) - \kme}_{\calH_k}  \leq  \calB_{k, \delta}(\bbP, n),
\end{equation*}
with
\begin{equation*}
    \calB_{k, \delta}(\bbP, n) = \calB_{ \delta} \eqdef \sqrt{2 v \frac{\log (2/ \delta)}{n}} + \frac{4 }{3}\sqrt{\overline{k}}\frac{\log (2/ \delta)}{n}.
\end{equation*}
\end{theorem}

\begin{remark}
\label{remark:superior-most-settings}
Since the reproducing kernel is bounded, it is always the case that
$$v \leq \bbE_{X \sim \bbP} \nrm{k(X, \cdot)}_{\calH_k}^2 \leq \bbE_{X \sim \bbP} k(X,X) \leq \overline{k},$$
where the first inequality can be found for example in 
\citep[Lemma~7.1, Proof]{cherief2022finite}.
As a result, Theorem~\ref{theorem:variance-aware-mmd-estimation} strictly supersedes \eqref{eq:tolstikhin-bound} at least when 
\begin{equation*}
n > \frac{16}{9} \left(\frac{\log (2 / \delta)}{1 + \sqrt{2 \log(1/\delta)} - \sqrt{2\log(2/\delta)}} \right)^2.
\end{equation*}
For instance, for $\delta = 0.05$, it suffices that $n \geq 46$.
\end{remark}

\begin{remark}
The Bernstein approach, in contrast to bounded differences, has the additional advantage of yielding a maximal inequality, further opening the door for early stopping methods.
\end{remark}

\begin{remark}[Sharpness of the constants]
We note but do not pursue that there exist other families of concentration inequalities which are known to dominate Bennett--Bernstein-type inequalities.
For instance, the class of inequalities pioneered by \citet{bentkus2006domination}, with empirical counterparts derived by \citet{kuchibhotla2021near}, which is known to be nearly optimal for sample averages of independent observations from a log-concave probability distribution. We also mention the family of inequalities obtained by \citet{waudby2024estimating} based on martingale methods.
However, our problem requires a bound for norms of sums of vectors in Hilbert spaces, or alternatively for the supremum of empirical processes, which we could not locate in the aforementioned references.
Additionally, even in the simpler case of sample averages, computation of the bound requires effort, thus a concentration bound in the more complicated Hilbert space setting could be challenging to compute. Finally, Bernstein-type bounds have known extensions for the case of time-dependent data, enabling us to analyze the estimation problem from stationary mixing sequences of observations (refer to Section~\ref{section:dependent-data}).
\end{remark}

We immediately observe that:
\begin{enumerate}[$(O1)$]
    \item While the bound in \eqref{eq:tolstikhin-bound} depends solely on chosen quantities and is computable without any knowledge of $\bbP$, this is not the case for $\calB_{\delta}$ in Theorem~\ref{theorem:variance-aware-mmd-estimation}. 
    \item Perhaps even more concerning, it is a priori unclear how $v$ depends on properties of $k$ and $\bbP$, thus making it difficult to convert assumptions on $\bbP$ into a bound for $v$.
\end{enumerate}
 We defer $(O1)$ to Section~\ref{section:variance-empirical} and first address $(O2)$ by pointing out that when $\calX = \R^d$, and for numerous hyper-parametrized families of kernels, it is possible to promptly obtain upper bounds on $v$ that decouple the influence of the hyper-parameter and some measure of spread of the underlying distribution. 

\subsection{Gaussian Kernel}
For a TI kernel---see \eqref{definition:translation-invariant}---we can rewrite $v$ as
\begin{equation*}
    v = \overline{\psi} - \E[X,X' \sim \bbP]{\psi(X' - X)}.
\end{equation*}
The Gaussian kernel with lengthscale parameter $\gamma > 0$, defined by
\begin{equation}
\label{definition:gaussian-kernel}
    k_\gamma(x,x') = \exp\left( - \frac{\nrm{x' - x}_2^2}{2 \gamma^2} \right),
\end{equation}
is the prototypical example of a characteristic translation invariant kernel,
and satisfies $\overline{\psi} = 1$. 
For $x \in \R^d$, let $x^i$ denote the $i$th component of $x$. The function $z \mapsto e^{-z}$ is convex, thus from Jensen's inequality 
$$v_\gamma = \bbE_{X \sim \bbP} \nrm{k_\gamma(X, \cdot) - \kme}_{\calH_{k_\gamma}}^2 \leq 1 - \exp \left( - \frac{\bbE_{X,X' \sim \bbP} \nrm{X' - X}_2^2}{2 \gamma^2} \right).$$
We can further rewrite the expectation on the right side of this inequality as
\begin{equation*}
\begin{split}
    \bbE_{X,X' \sim \bbP} \nrm{X' - X}_2^2 &=\sum_{i = 1}^{d} \E[X,X' \sim \bbP]{  (X'^i-X^i)^2} = \sum_{i = 1}^{d} 2\bbV_{X \sim \bbP}  X^i = 2 \Trace \Sigma_{\bbP}, \
\end{split}
\end{equation*}
where $\Sigma_{\bbP} \eqdef \bbV_{X \sim \bbP} X = \E[X \sim \bbP]{(X - \bbE X)(X - \bbE X)^\star}$ is the covariance matrix of $\bbP$ and its trace is interpreted as a measure a total variance, agnostic to correlations between distinct components.
As a result, for any fixed $\gamma > 0$, we obtain from Theorem~\ref{theorem:variance-aware-mmd-estimation} that with probability at least $1 -\delta$,
\begin{equation}
\label{eq:variance-aware-gaussian-mmd}
    \nrm{\ekme - \kme}_{\calH_{k_\gamma}} \leq \sqrt{ 2 \left(1 - e^{- \Trace \Sigma_{\bbP} / \gamma^2} \right) \frac{\log (2/ \delta)}{n}} + \frac{4}{3} \frac{\log (2/ \delta)}{n}.
\end{equation}
While $\gamma$ has no influence over the right-hand side in \eqref{eq:tolstikhin-bound}, it is clear that for $\gamma \to \infty$, the rate of convergence will be accelerated in \eqref{eq:variance-aware-gaussian-mmd}.

\begin{example}[Gaussian location model with known variance]
\label{example:gaussian-known-variance}
Assume that $\bbP = \bbP_\theta = \calN(\theta, \sigma^2 I_d)$, the Gaussian distribution with unknown location parameter $\theta \in \R^d$, but with known covariance matrix $\sigma^2 I_d$.
It holds that $\Trace \Sigma_{\bbP_\theta} = \sigma^2 d$. 
\begin{enumerate}[$(i)$]
    \item Fixing $\gamma$ and taking $\sigma \to 0$, the bound in \eqref{eq:variance-aware-gaussian-mmd} vanishes, unlike \eqref{eq:tolstikhin-bound}.
    \item Setting $\gamma^2 = \lambda 
\sigma^2 d$ with $\lambda > 0$, we readily obtain the variance upper bound $v \leq (1 - e^{-1/\lambda}) \leq 1/\lambda$, enabling uncomplicated tuning of the convergence rate by $\lambda$.
\end{enumerate}
\end{example}

\begin{remark}
\label{remark:arbitrary-rate}
Example~\ref{example:gaussian-known-variance} highlights that if we allow the lengthscale parameter to vary with the sample size $n$,
Theorem~\ref{theorem:variance-aware-mmd-estimation} speeds up the convergence rate dramatically. For instance, setting $\lambda_n = -1 / \log (1 - 1/n)$, we achieve a rate in $\bigO_P(n^{-1})$, subverting the lower bounds of \citet{tolstikhin2017minimax}.
In fact, we can reach any prescribed rate of convergence between $\bigO_P(n^{-1/2})$ and $\bigO_P(n^{-1})$ for a suitable choice of $\lambda_n$. This stands in sharp contrast with bounds obtained through a bounded differences approach such as \eqref{eq:tolstikhin-bound}. However, it is important to remember that a larger lengthscale parameter will flatten the kernel, hence make the left hand side in \eqref{eq:variance-aware-gaussian-mmd} less informative. Depending on the application, it is possible to achieve the optimal balance between these two effects; some examples are discussed in Section~\ref{section:applications}.
\end{remark}

\subsection{Convex Radial Square Basis Functions}

In fact, a technique similar to that used for obtaining \eqref{eq:variance-aware-gaussian-mmd} yields a bound on the variance for any (convex) radial kernel.

\begin{lemma}
Assume that for any $x,x' \in \calX$, $k(x,x') = r(-\nrm{x' - x}_2^2)$ for some convex function $r \colon \R_+ \to \R$. Then
\begin{equation*}
    v \leq \overline{r} - r\left(-\Trace  \Sigma_{\bbP} \right),
\end{equation*}
where $\Sigma_{\bbP}$ is the covariance matrix that pertains to $\bbP$.
\end{lemma}

\subsection{Positive Definitive Matrix on the Finite Space}

Consider $\abs{\calX} < \infty$ and a symmetric positive definite matrix $K$ of size $\abs{\calX}$.
Write $K = \sum_{x \in \calX} \lambda_x v_x \trn v_x$, where for any $x \in \calX$, $K v_x = \lambda_x v_x$, and by positive definiteness, $\lambda_x >0$.
Then the feature map can be expressed as
\begin{equation*}
    K(x, \cdot) = \left( \sqrt{\lambda_{x'}} v_x(x') \right)_{x' \in \calX},
\end{equation*}
and by direct calculation,
\begin{equation*}
\begin{split}
    \nrm{\ekme - \kme}_{\calH_K}^2 &= \sum_{x \in \calX} \lambda_x \left( \frac{1}{n} \sum_{t=1}^{n} v_x(X_t) - \bbE  [v_x(X)] \right)^2, \\
    v &= \sum_{x \in \calX} \lambda_x \bbV [v_x(X)].
    \end{split}
\end{equation*}
In particular, for $K = I$,
\begin{equation*}
\begin{split}
    \nrm{\ekme - \kme}_{\calH_K} = \nrm{\widehat{\bbP}_n - \bbP}_2, \qquad v = 1 - \nrm{\bbP}_2^2,
    \end{split}
\end{equation*}
recovering that with probability at least $1- \delta$,
\begin{equation*}
    \nrm{\widehat{\bbP}_n - \bbP}_2 \leq \sqrt{2(1 - \nrm{\bbP}_2^2)\frac{\log (2/\delta)}{n}} + \frac{4}{3} \frac{\log (2/\delta)}{n}.
\end{equation*}

So far, we have shown how to obtain for a natural class of kernels an upper bound on $v$ that decouples the choice of $k$ from some measure of total variance of $\bbP$.
Provided an upper bound on the latter is available (see Example~\ref{example:gaussian-known-variance}), we obtain an improved convergence rate.
However, 
we understandably may not have a bound on the variance of $\bbP$, or could have insight about its variance, but only have access to contaminated data (see later Section~\ref{section:improving-bounds-contaminated-data}).
Recovering a bound that still enjoys the discussed speed-up without a priori knowledge on $\bbP$ is the problem we set out to solve in the next section.

\section{Convergence Rates with Empirical Variance Proxy}
\label{section:variance-empirical}
We solve the issue $(O1)$ of not knowing $v$ by estimating it from the data, simultaneously to $\mu_\bbP$.
A short heuristic analysis of the 
Epanechnikov function (which is not typically a kernel according to our definition) that we conduct in Section~\ref{section:variance-empirical-intuition} hints at an ``empirical Bernstein" \citep{audibert2007tuning, maurer2009empirical} approach---replacing the variance term by some empirical proxy---which we explore formally in Section~\ref{section:variance-empirical-systematic}.
 The argument is structured 
 around
the pivotal definition
of a weakly self-bounding function that we first recall.

\begin{definition}[{\citealt{maurer2006concentration, boucheron2009concentration}}]
\label{definition:weakly-self-bounding-and-bounded-differences-function}

Let $n \in \N$, $t_0 \in [n]$,  
$$x = (x_1, \dots, x_{t_0}, \dots, x_n) \in \calX^n,$$ $x_{t_0}' \in \calX$, and write 
$$x^{(t_0)} = (x_1, \dots, x_{{t_0}-1}, x_{t_0}', x_{{t_0}+1}, \dots, x_n),$$ for the vector where $x_{t_0}$ has been replaced with $x_{t_0'}$.
Let $f \colon \calX^n \to \R$.

\begin{enumerate}[$(i)$]
    \item 
The function $f$ is called weakly $(\alpha, \beta)$-self-bounding
when for all $x \in \calX^n$,
\begin{equation*}
    \sum_{t_0 = 1}^{n} \left( f(x) - \inf_{x_{t_0}' \in \calX} f\left(x^{(t_0)}\right)\right)^2 \leq \alpha f(x) + \beta.
\end{equation*}
\item The function $f$ is said to have the bounded differences property when for all $x \in \calX^n$ and all $t_0 \in [n]$,
\begin{equation*}
    f(x) - \inf_{x_{t_0}' \in \calX} f\left(x^{(t_0)}\right) \leq 1.
\end{equation*}
\end{enumerate}
\end{definition}

\subsection{Intuition in the Hypercube}
\label{section:variance-empirical-intuition}
Let us take $\calX = [0,1]^d$ the binary hypercube for $d \gg 1$, and consider
the Epanechnikov (parabolic) function
\begin{equation*}
    q(x,x') =  1 -  \nrm{x' - x}_2^2 / d.
\end{equation*}
Note that $q$ does not define a proper kernel \citep[p.9]{cuturi2009positive}.
Nevertheless, the function $q$ is TI, with $\psi(t) = 1 - \nrm{t}_2^2/d$, $\overline{\psi} = 1$, and we can extend the definition of the variance term
\begin{equation*}
\begin{split}
    v &= \overline{\psi} - \E[X,X' \sim \bbP]{\psi(X'-X)} = 2 \Trace  \Sigma_\bbP /d.
\end{split}
\end{equation*} 
It is natural to define for $i \in [d]$,
\begin{equation*}
\widehat{v}^i(X_1, \dots, X_n) \eqdef \frac{1}{2 n (n-1)} \sum_{s = 1}^{n}\sum_{t = 1}^{n} \left(X_t^i - X_s^i\right)^2,
\end{equation*}
that will act as an unbiased empirical proxy for $v^i \eqdef \Var[X \sim \bbP]{X^i}$, and introduce
$\Trace  \widehat{\Sigma} \eqdef \sum_{i =1}^{d} \widehat{v}^i$,
as an estimator for the trace of the covariance matrix $\Sigma_{\bbP}$. 

\begin{lemma}
\label{lemma:epanechnikov-self-bounding}
For $x \in [0,1]^d$, the function $\calX^m \to \R_+, x \mapsto \Trace  \widehat{\Sigma}(x) /d$
is weakly $(n/(n-1), 0)$-self-bounding and has the bounded differences property in the sense of Definition~\ref{definition:weakly-self-bounding-and-bounded-differences-function}.
\end{lemma}
As a consequence of Lemma~\ref{lemma:epanechnikov-self-bounding}, the technique of  \citet[Theorem~10]{maurer2009empirical} provides the following deviation bounds on the square root of the total variance
\begin{equation}
\label{eq:total-variance-deviation-bounds}
\begin{split}
\PR{ b \left[ \sqrt{\Trace  \Sigma_{\bbP}/d} - \sqrt{\Trace  \widehat{\Sigma}(X)/d} \right] > \sqrt{\frac{2 \log (1/ \delta)}{n - 1}} } \leq \delta, \qquad b \in \set{-1, 1}.
\end{split}
\end{equation}
In other words, with high confidence, we can replace the trace of the covariance matrix with its empirical proxy in the convergence bounds.

\subsection{Systematic Approach}
\label{section:variance-empirical-systematic}
Our goal is to rigorously  develop the approach intuited in the previous section to include a large class of reproducing kernels.
We propose the following empirical proxy for the variance term $v$,
\begin{equation}
\label{eq:general-empirical-variance-proxy}
    \widehat{v}_k(X_1, \dots, X_n) \eqdef \frac{1}{n-1} \sum_{t = 1}^{n} \left( k(X_t, X_t) - \frac{1}{n}\sum_{s=1}^{n} k(X_t, X_s) \right),
\end{equation}
write more simply $\widehat{v}$ for $\widehat{v}_k$ when $k$ is clear from the context,
and promptly verify that $\widehat{v}$ is an unbiased estimator for $v$.
\begin{lemma}[Unbiasedness]
\label{lemma:general-empirical-proxy-unbiased}
It holds that 
$\bbE \widehat{v} = v$,
where the expectation is taken over the sample $X_1, \dots, X_n \stackrel{iid}{\sim} \bbP$.
\end{lemma}

The remainder of this section is devoted to analyzing the self-boundedness properties of $\widehat{v}$ under mild conditions on $k$ and deriving corresponding empirical confidence intervals around $\kme$.
We henceforth let $k$ be a characteristic TI kernel defined by the positive definite function $\psi$.
We also obtain partial results for the more general case of kernels which are not translation invariant. This analysis is deferred to the appendix Section~\ref{section:non-translation-invariant-kernels} (refer to Theorem~\ref{theorem:general-empirical-variance-mmd-bound}).
In the TI setting, the expression of $\widehat{v}$ can be simplified as
\begin{equation}
\label{eq:translation-invariant-empirical-variance-proxy}
    \widehat{v}(X_1, \dots, X_n) \eqdef \overline{\psi} - \frac{1}{(n-1)n} \sum_{t \neq s} \psi (X_t - X_s).
\end{equation}
Since the kernel is characteristic, $\psi$ cannot be a constant function, that is $\Delta \psi > 0$.
Under our assumptions, we introduce a function involving $\widehat{v}$ and $\Delta \psi$ that is weakly self-bounded and has the bounded differences property.
\begin{lemma}
\label{lemma:translation-invariant-weakly-self-bounding}
Let $\psi$ define a characteristic TI kernel.
The function
\begin{equation*}
    \calX^n \to \R, \, x \mapsto \frac{n}{2 \Delta \psi} \widehat{v}(x),
\end{equation*}
is weakly $(2,0)$-self-bounding and has bounded differences in the sense of Definition~\ref{definition:weakly-self-bounding-and-bounded-differences-function}.
\end{lemma}
This property leads to concentration of $\sqrt{\widehat{v}}$ around $\sqrt{v}$.
\begin{lemma}
\label{lemma:concentration-square-root-v}
For $b \in \set{-1,1}$, with probability at least $1- \delta$,
\begin{equation}
\label{eq:square-root-confidence-interval-for-v}
    b \left[ \sqrt{\widehat{v}} - \sqrt{v} \right] \leq 2 \sqrt{\frac{2 \Delta \psi \log (1 /\delta)}{n}}.
\end{equation}
\end{lemma}

\begin{theorem}[Confidence interval with empirical variance for TI kernel]

\label{theorem:empirical-variance-mmd-bound}
Let $n \in \N$, $X_1, \dots, X_n \stackrel{iid}{\sim} \bbP$, let $k$ be a characteristic translation invariant kernel defined from a positive definite function $\psi$ [see \eqref{definition:translation-invariant}].
Then with probability at least $1- \delta$, it holds that
\begin{equation*}
\begin{split}
    \nrm{\ekme - \kme}_{\calH_{k}} &\leq \widehat{\calB}_{k, \delta}(X_1, \dots, X_n), \\
    \end{split}
    \end{equation*}
    with 
    \begin{equation*}
\begin{split}
    \widehat{\calB}_{k, \delta}(X_1, \dots, X_n) = \widehat{\calB}_{\delta} \eqdef \sqrt{2 \widehat{v}(X_1,\dots, X_n) \frac{\log (4/ \delta)}{n}} + \frac{16}{3} \sqrt{\Delta \psi} \frac{\log (4/ \delta)}{n},
\end{split}
\end{equation*}
and where $\widehat{v}$ is the empirical variance proxy defined in \eqref{eq:translation-invariant-empirical-variance-proxy}.
\end{theorem}

\begin{example}
Theorem~\ref{theorem:empirical-variance-mmd-bound} immediately holds for Gaussian kernels, regardless of the lengthscale parameter.
\end{example}

\begin{remark}
\label{remark:observations-on-v}
We make the following observations.
\begin{enumerate}[$(i)$]
    \item It still almost surely holds that $\widehat{v} \leq \Delta k \leq 2\overline{k}$, thus the fully empirical bound can never be more than a constant factor away from \eqref{eq:tolstikhin-bound} or Theorem~\ref{theorem:variance-aware-mmd-estimation}.
    \item Crucially, self-boundedness of $\widehat{v}$ and being able to apply Theorem~\ref{theorem:empirical-variance-mmd-bound} (or Theorem~\ref{theorem:general-empirical-variance-mmd-bound} for non TI kernels) only depends on the choice of kernel, and not on the properties of the underlying distribution.
    \item Theorem~\ref{theorem:empirical-variance-mmd-bound} does not depend on smoothness properties of the kernel (also true for Theorem~\ref{theorem:general-empirical-variance-mmd-bound} for non TI kernels).
\end{enumerate}
\end{remark}

\subsection{Computability of the Empirical Variance Proxy}

The proxy $\widehat{v}$ is computable from data with $\bigO(n^2)$ calls to the kernel function. In Section~\ref{section:applications}, we will discuss how to use variance-aware bounds to improve confidence bounds and test procedures based on the minimum-MMD estimator. In general, this estimator is computed by stochastic gradient descent (SGD)~\citep{dziugaite2015training,li2015generative,cherief2022finite} or variants like stochastic natural gradient descent~\citep{briol2019statistical}. A single step of SGD requires to sample $m$ iid random variables from $\bbP_{\theta}$ where $\theta$ is the current estimate, and to compute the unbiased estimate of the gradient which requires (among others) $\bigO(mn)$ calls to the kernel function.  In the natural version, one must also compute a Jacobian at each step, which requires $\bigO(n^2)$ calls to the kernel. Morever, in the discussion following Theorem 2 of \cite{briol2019statistical}, it is argued that we should take $m\gtrsim n$ in these algorithms. Thus, both SGD and its natural variant will require at least $\bigO(n^2)$ calls to the kernel at each iteration. In this light, the computation of $\widehat{v}$ does not significantly affect the computational burden. Note however that, in a few situations where the gradient of the MMD is available in close form, it is possible to use a non-stochastic gradient descent. Such examples include Gaussian mean estimation~\citep{cherief2022finite} and Gaussian copula estimation~\citep{alquier2022estimation}. Each step of the gradient descent requires only $\bigO(n)$ calls to the kernel and convergence is typically very fast. In these situations, the computation of $\widehat{v}$ can  increase the computational cost when $n$ is large.

\section{Convergence Rates for the Difference of Two Means}
\label{section:two-means}
We now extend the results of the previous section to the estimation of the norm of the difference of two means $\nrm{\mu_{\bbP}-\mu_{\bbQ}}_{\calH_k}$. One of the applications of these results is two-sample tests, where we test the hypothesis $\bbP=\bbQ$, as in \cite{gretton2012kernel}. We will cover this application in Section~\ref{section:applications}. However, they also have an interest on their own, to illustrate the benefits of variance-aware bounds.

In all this section, we assume that two samples are observed:  $X_1,\dots,X_n$ iid from $\bbP$, and $Y_1,\dots,Y_m$ iid from $\bbQ$. We want to estimate the norm of the difference $\nrm{\mu_{\bbP}-\mu_{\bbQ}}_{\calH_k}$.

\subsection{First Approach: Estimation by a U-Statistic}

A classical approach is to estimate the squared norm $\nrm{\mu_{\bbP}-\mu_{\bbQ}}_{\calH_k}^2$ using a U-statistic, as proposed by~\cite{gretton2012kernel}, and then to use a Hoeffding or Bernstein-type inequality for U-statistics. 

While~\cite{gretton2012kernel} used the U-statistics version of Hoeffding's inequality to control the risk of this estimator, using a Bernstein inequality instead might allow to derive variance-aware bounds. The first version of such inequalities were proven by~\cite{hoeffding1963probability} and~\cite{arcones1995bernstein}. More recently, the results of these two papers are included in Theorem 2 of~\cite{peel2010empirical} while Theorem 3 gives an empirical version of each.

We now provide more details on this approach. First, note that
$$ 
\nrm{\mu_{\bbP}-\mu_{\bbQ}}_{\calH_k}^2 = \bbE_{X,X' \sim \bbP} k(X,X') + \bbE_{Y,Y' \sim \bbQ} k(Y,Y')-2 \bbE_{X \sim \bbP,Y\sim\bbQ} k(X,Y)
$$
can be directly estimated by the U-statistic:
$$
\tilde{U}_{n,m} = \frac{1}{n(n-1)} \sum_{i=1}^n \sum_{j\neq i} k(X_i,X_j) +  \frac{1}{m(m-1)} \sum_{i=1}^m \sum_{j\neq i} k(Y_i,Y_j) - \frac{2}{mn} \sum_{i=1}^n \sum_{j=1}^m k(X_i,Y_j)
$$
as proposed by~\cite{gretton2012kernel}. This is an unbiased estimator: $\mathbb{E}\tilde{U}_{n,m} = \nrm{\mu_{\bbP}-\mu_{\bbQ}}_{\calH_k}^2$. In the special case $m=n$, we can replace $\tilde{U}_{n,n}$ by the simpler
$$
\hat{U}_n = \frac{1}{n(n-1)}  \sum_{i=1}^n \sum_{j\neq i} \left[k(X_i,X_j) + k(Y_i,Y_j) - k(X_i,Y_j) - k (X_j,Y_i) \right].
$$
We will illustrate the U-statistics approach in this simpler setting.

Using Theorem 2 in~\cite{peel2010empirical}, with their $q((X_i,Y_i),(X_j,Y_j)) = k(X_i,X_j) + k(Y_i,Y_j) - k(X_i,Y_j) - k (X_j,Y_i)$ and their $m=2$, we obtain the first part of the following statement. Using Theorem 3, we obtain the second one.
\begin{theorem}
\label{thm:u-statistics}
Put
$$\sigma^2_{\bbP,\bbQ} = \mathbb{V}_{X\sim\bbP,Y\sim \bbQ}[ \mathbb{E}_{X'\sim\bbP} k(X,X') +  \mathbb{E}_{Y'\sim\bbQ} k(Y,Y') - \mathbb{E}_{Y'\sim\bbQ} k(X,Y') - \mathbb{E}_{X'\sim\bbP} k(X',Y)  ]  $$
and
\begin{align*}
\hat{\sigma}^2_{\bbP,\bbQ} = \frac{1}{n(n-1)(n-2)}
\sum_{i=1}^n \sum_{j\neq i} \sum_{k\neq i,j}
&
\biggl[k(X_i,X_j) + k(Y_i,Y_j) - k(X_i,Y_j) - k (X_j,Y_i) \biggr]
\\
&
\times \biggl[k(X_i,X_k) + k(Y_i,Y_k) - k(X_i,Y_k) - k (X_k,Y_i) \biggr]
.
\end{align*}
Note that $\mathbb{E}\hat{\sigma}^2_{\bbP,\bbQ} = \sigma^2_{\bbP,\bbQ} + \nrm{\mu_{\bbP}-\mu_{\bbQ}}_{\calH_k}^4 $.
For any $\delta\in(0,1]$, with probability at least $1-\delta$,
$$
\left| \nrm{\mu_{\bbP}-\mu_{\bbQ}}_{\calH_k}^2 - \hat{U}_n\right|
\leq
2\Delta k \sqrt{ \frac{8 \sigma^2_{\bbP,\bbQ} }{n} \log\frac{4}{\delta} } + 2\Delta k  \frac{64 + \frac{1}{6}}{n}\log\frac{4}{\delta}.
$$
Moreover, with probability at least $1-\delta$,
$$
\left| \nrm{\mu_{\bbP}-\mu_{\bbQ}}_{\calH_k}^2 - \hat{U}_n\right|
\leq
2\Delta k  \sqrt{ \frac{8 \hat{\sigma}^2_{\bbP,\bbQ} }{n} \log\frac{8}{\delta} } + 2\Delta k \frac{64 + \frac{1}{6} + 5\sqrt{2}}{n}\log\frac{8}{\delta}.
$$
\end{theorem}
Observe that, when $\bbP\neq\bbQ$, the bound states that $\sqrt{\hat{U}_n}$ will be of the order of $ \nrm{\mu_{\bbP}-\mu_{\bbQ}}_{\calH_k}$. On the other hand, when $\bbP=\bbQ$, we have both $ \nrm{\mu_{\bbP}-\mu_{\bbQ}}_{\calH_k}=0$ and $\sigma^2_{\bbP,\bbQ} = 0$, and thus the first inequality in the theorem gives:
\begin{equation}
\label{eq:ustat1}
 \sqrt{\hat{U}_n} \leq \sqrt{ 2\Delta k  \frac{64 + \frac{1}{6}+5\sqrt{2}}{n}\log\frac{8}{\delta} },
\end{equation}
which does not depend on the variances of $\bbP$ and $\bbQ$.

\subsection{Variance-Aware Control of the Fluctuations for Each Sample}

An alternative approach is to apply the triangle inequality to upper bound separately the fluctuations for each sample:
\begin{multline*}
\biggl|
\nrm{\mu_{\bbP}-\mu_{\bbQ}}_{\calH_k}
-
\nrm{\hat{\mu}_{\bbP}(X_1,\dots,X_n)-\hat{\mu}_{\bbQ}(Y_1,\dots,Y_n) }_{\calH_k}
\biggr|
\\
\leq \nrm{\mu_{\bbP}-\hat{\mu}_{\bbP}(X_1,\dots,X_n)}_{\calH_k}
+
\nrm{\mu_{\bbQ}-\hat{\mu}_{\bbQ}(Y_1,\dots,Y_n) }_{\calH_k}.
\end{multline*}
For example, a direct application of Theorem~\ref{theorem:variance-aware-mmd-estimation} with a union bound gives the following corollary.
\begin{corollary}
\label{coro:u-statistics}
Let $X_1,\dots,X_n\sim\bbP$, $Y_1,\dots,Y_m\sim \bbQ$, $k \colon \calX \times \calX \to \R$ be a reproducing kernel, and  let
\begin{equation*}
    v_{\bbP} \eqdef \bbE_{X \sim \bbP} \nrm{k(X, \cdot) - \mu_{\bbP}}_{\calH_k}^2 \text{ and }  v_{\bbQ} \eqdef \bbE_{Y \sim \bbQ} \nrm{k(Y, \cdot) - \mu_{\bbQ}}_{\calH_k}^2
\end{equation*}
and $\overline{k} = \sup_{x \in\calX} k(x,x) $. With probability at least $1 - \delta$, it holds that
\begin{multline*}
  \biggl|
\nrm{\mu_{\bbP}-\mu_{\bbQ}}_{\calH_k}
-
\nrm{\hat{\mu}_{\bbP}(X_1,\dots,X_n)-\hat{\mu}_{\bbQ}(Y_1,\dots,Y_m) }_{\calH_k}
\biggr| 
\\
\leq
\sqrt{2 v_\bbP \frac{\log (4/ \delta)}{n}} + \sqrt{2 v_\bbQ \frac{\log (4/ \delta)}{m}} + \left(\frac{1}{n} + \frac{1}{m}\right)\frac{4 }{3}\sqrt{\overline{k}}\log (4/ \delta).
\end{multline*}
\end{corollary}
Of course, we can also state results with empirical variance instead of the true variance, by using Theorems~\ref{theorem:empirical-variance-mmd-bound} and~\ref{theorem:general-empirical-variance-mmd-bound}.

In order to compare this result to the U-statistics approach, consider the case $n=m$. Corollary~\ref{coro:u-statistics} gives the following upper bound:
\begin{multline*}
\nrm{\hat{\mu}_{\bbP}(X_1,\dots,X_n)-\hat{\mu}_{\bbQ}(Y_1,\dots,Y_n) }_{\calH_k}
\\
\leq 
\nrm{\mu_{\bbP}-\mu_{\bbQ}}_{\calH_k}
+
\sqrt{2 v_\bbP \frac{\log (4/ \delta)}{n}} + \sqrt{2 v_\bbQ \frac{\log (4/ \delta)}{n}} + \frac{8 }{3n}\sqrt{\overline{k}}\log (4/ \delta).
\end{multline*}
In particular, when $\bbP=\bbQ$, the bound becomes:
\begin{equation}
\label{eq:ustat2}
\nrm{\hat{\mu}_{\bbP}(X_1,\dots,X_n)-\hat{\mu}_{\bbQ}(Y_1,\dots,Y_n) }_{\calH_k}
\leq 2 \sqrt{2 v_\bbP \frac{\log (4/ \delta)}{n}} + \frac{8 }{3n}\sqrt{\overline{k}}\log (2/ \delta).
\end{equation}
Even though $\nrm{\hat{\mu}_{\bbP}(X_1,\dots,X_n)-\hat{\mu}_{\bbQ}(Y_1,\dots,Y_n) }_{\calH_k}^2$ is not an unbiased estimator of $\nrm{\mu_{\bbP}-\mu_{\bbQ}}_{\calH_k}^2$ as $\hat{U}$,~\eqref{eq:ustat1} and~\eqref{eq:ustat2} show the fluctuations of $\nrm{\hat{\mu}_{\bbP}(X_1,\dots,X_n)-\hat{\mu}_{\bbQ}(Y_1,\dots,Y_n) }_{\calH_k}$ when $\bbP=\bbQ$ are upper bounded by $\sqrt{v_{\bbP}/n}+1/n $ while the ones of $\sqrt{\hat{U}}$ by $\sqrt{1/n}$. This can be a serious improvement if $v_\bbP$ is small. Note that we do not claim superiority of the plug-in estimator, but rather that the currently available non-asymptotic bounds for this estimator are tighter. This is an illustration of the power of the variance-aware bounds.
A way to compare both estimators would be through an accurate study of their asymptotic fluctuations when $\bbP=\bbQ$. So far, this analysis is only available for the U-statistic, see Theorem~12 in~\citet{gretton2012kernel}.

\section{Convergence Rates with Time-Dependent Data}
\label{section:dependent-data}
In this section, we establish convergence rates for cases where the data $X_1, \dots, X_n $ is not independent.
Namely, we will assume the data to be a stationary mixing sequence \citep{bradley2005basic, doukhan2012mixing}. In this setting, the observations are identically distributed with marginal distribution $\bbP$, but  exhibit time dependencies that diminish as the time interval increases.
\begin{lemma}
\label{lemma:exact-expression-second-moment-covariances}
Suppose that $X_1, \dots, X_n$ is a stationary sequence of possibly dependent random variables.
Then
\begin{equation*}
\begin{split}
    \bbE\nrm{ \ekme - \kme }^2_{\calH_k}
    &= \frac{1}{n}\left(v + \Sigma_n \right), \\
\end{split}
\end{equation*}
with
\begin{equation*}
    \Sigma_n \eqdef \frac{2}{n}\sum_{t = 2}^{n} (n - t + 1) \rho_t, \qquad \text{ where } \qquad \rho_t \eqdef \bbE \langle{ k(X_t, \cdot) - \kme,  k(X_1, \cdot) - \kme \rangle}_{\calH_k},
\end{equation*}
are the covariance coefficients in the RKHS introduced in \citet{cherief2022finite}. 
In particular, when the process is iid, the expression above simplifies to $$\bbE\nrm{ \ekme - \kme }^2_{\calH_k} = \frac{v}{n}.$$
\end{lemma}

\subsection{For \texorpdfstring{$\phi$}{Phi}-Mixing Processes}

The first flavor of mixing we consider is $\phi$-mixing, introduced by \citet{ibragimov1962some}.
Recall that the $\phi$-mixing coefficient
\citep{bradley2005basic, doukhan2012mixing} is defined for two $\sigma$-fields $\calA$ and $\calB$ by
 \begin{equation*}
     \phi(\calA, \calB) \eqdef \sup_{\substack{A \in \calA, B \in \calB \\ \PR{A} > 0}} \abs{\PR{B | A} - \PR{B}}.
 \end{equation*}
For $s \in \bbN$, we further define
\begin{equation*}
     \phi(s) \eqdef \sup_{r \in \bbN} \phi \left( \sigma\left( \set{X_t \colon t \leq r} \right), \sigma\left( \set{X_t \colon t \geq r + s} \right) \right),
 \end{equation*}
 where for $T \subset \bbN$, $\sigma(\{X_t\}_{t \in T})$ is the $\sigma$-field generated by the random variables $\{X_t\}_{t \in T}$.
The random process is then called $\phi$-mixing when $\lim_{s \to \infty} \phi(s) = 0$.
Additionally, we define the triangular coupling matrix $\Gamma$ \citep{sampson2000concentration} as follows. For $t, s \in [n]$,
\begin{equation}
\label{definition:coupling-matrix}
    \Gamma(t,s) \eqdef \begin{cases}
        1 &\text{ when } t = s \\
        0 &\text{ when } t > s \\
        \sqrt{2 \phi(\sigma(X_1, \dots, X_t), \sigma(X_s, \dots, X_n))} &\text{ otherwise}.
    \end{cases}
\end{equation} 

\begin{theorem}[Variance-aware confidence interval with $\phi$-mixing data]
\label{theorem:variance-aware-mmd-estimation-phi-mixing}
    Let $\delta \in (0,1)$, and $X_1, \dots, X_n$ be a stationary $\phi$-mixing sequence with marginal distribution $\bbP$. We let $k \colon \calX \times \calX \to \R$ be a reproducing kernel with $\overline{k} = \sup_{x \in\calX} k(x,x)$. 
With probability at least $1 - \delta$, it holds that
\begin{equation*}
    \nrm{\ekme(X_1, \dots, X_n) - \kme}_{\calH_k}  \leq  \calB^{\phi}_{k, \delta}(\bbP, n),
\end{equation*}
with
\begin{equation*}
    \calB^{\phi}_{k, \delta}(\bbP, n) = \calB^{\phi}_{\delta} \eqdef \sqrt{\frac{v + \Sigma_n}{n}} + 4 \sqrt{\frac{2 v \nrm{\Gamma}_2 \log (1 /\delta)}{n}} + \frac{8 \overline{k} \nrm{\Gamma}_2 \log (1/\delta)}{n},
\end{equation*}
where $\Sigma_n$ is defined in Lemma~\ref{lemma:exact-expression-second-moment-covariances}, $\Gamma$ is defined in Eq.~\eqref{definition:coupling-matrix} and $\nrm{\cdot}_2$ is the spectral norm.
\end{theorem}

In particular, when the process is iid, $\Gamma$ is the identity matrix, hence $\nrm{\Gamma}_2 = 1$.

\begin{example}[Uniformly ergodic Markov chains]
Suppose that there exists $\phi_0, \phi_1 \in \bbR_+$ such that for all $t \in \bbR$, $\phi(t) \leq \phi_0 \exp( - \phi_1 t)$. Then it holds (refer for example to \citet{sampson2000concentration}) that
\begin{equation*}
    \nrm{\Gamma}_2 \leq \frac{\sqrt{2 \phi_0}}{ 1 - e^{-\phi_1 / 2}}.
\end{equation*}
In particular, this is known to hold for a uniformly ergodic Markov chain \citep{10.1214/154957804100000024} with transition operator $P$ and stationary distribution $\pi$. In this case we can choose, $\phi_0 = 4$ and $\phi_1 = \log(2)/\tmix$, where 
\begin{equation*}
    \tmix \eqdef \min_{t \in \bbN}
    \set{\sup_{x \in \calX} \tv{ P^t(x , \cdot) - \pi(\cdot)} < 1/4},
\end{equation*}
is called the mixing time of the chain with respect to total variation \citep{levin2009markov}; here we obtain $\nrm{\Gamma}_2 \leq 10 \tmix$.  
\end{example}

\subsection{For \texorpdfstring{$\beta$}{Beta}-Mixing Processes}

We next consider $\beta$-mixing, a common assumption in the machine learning literature \citep{mohri2010stability}.
Recall that the $\beta$-mixing coefficient
\citep{bradley2005basic, doukhan2012mixing} is defined for two $\sigma$-fields $\calA$ and $\calB$,
 \begin{equation*}
     \beta(\calA, \calB) \eqdef \sup \frac{1}{2} \sum_{i = 1}^{I} \sum_{j = 1}^{J} \abs{\bbP(A_i \cap B_j) - \bbP(A_i)\bbP(B_j)},
 \end{equation*}
where the supremum is taken over all pairs of finite partitions $\{ A_1, \dots, A_I \}$
and $\{ B_1, \dots, B_J \}$ of $\calX$ such that for any $1 \leq i \leq I$ and any $1 \leq j \leq J$, $A_i \in \calA$  and $B_j \in \calB$ \citep[Equation 7]{bradley2005basic}.
For $s \in \bbN$, we further define
\begin{equation*}
     \beta(s) \eqdef \sup_{r \in \bbN} \beta \left( \sigma\left( \set{X_t \colon t \leq r} \right), \sigma\left( \set{X_t \colon t \geq r + s} \right) \right),
 \end{equation*}
 where for $T \subset \bbN$, $\sigma(\{X_t\}_{t \in T})$ is the $\sigma$-field generated by the random variables $\{X_t\}_{t \in T}$.
The random process is then called $\beta$-mixing\footnote{A $\beta$-mixing process is also called \emph{absolutely regular} in the literature.} when $\lim_{s \to \infty} \beta(s) = 0$.
In this case, for $\xi \in (0, 1)$, we define the $\beta$-mixing time as
\begin{equation*}
    \betamix(\xi) \eqdef \argmin_{t \in \bbN} \set{ \beta(t) < \xi }.
\end{equation*}
We observe that when the process is a stationary time-homogeneous Markov chain, $\betamix$ corresponds to the definition of the average-mixing time \citep{munch2023mixing, wolfer2024optimistic}.
Note that $\beta(\calA, \calB) \leq \phi(\calA, \calB)$ \citep[Equation 1.11]{bradley2005basic},
and no reverse inequality holds for any universal constant,
thus $\beta$-mixing is a strictly weaker notion than $\phi$-mixing.

\begin{theorem}[Variance-aware confidence interval with $\beta$-mixing data]
\label{theorem:variance-aware-mmd-estimation-beta-mixing}
    We let $\delta \in (0,1)$, and $X_1, \dots, X_n$ be a stationary $\beta$-mixing sequence. We let $k \colon \calX \times \calX \to \R$ be a reproducing kernel with $\overline{k} = \sup_{x \in\calX} k(x,x)$. 
    Finally, we suppose for simplicity\footnote{An adaptation of the proof removes this assumption and recovers similar bounds up to universal constants.} that $n$ is a multiple of $\tau^{\beta}_{n, \delta}$ where
    \begin{equation*}
    \tau^{\beta}_{n, \delta} \eqdef \argmin_{s \in \bbN} \set{ s \geq \betamix \left( \frac{\delta}{6(n/(2s) - 1)} \right)}.
\end{equation*}
With probability at least $1 - \delta$, it holds that
\begin{equation*}
    \nrm{\ekme(X_1, \dots, X_n) - \kme}_{\calH_k}  \leq  \calB^{\beta}_{k, \delta}(\bbP, n),
\end{equation*}
with
\begin{equation*}
    \calB^{\beta}_{k, \delta}(\bbP, n) = \calB^{\beta}_{ \delta} \eqdef 2 \sqrt{ \left(v + \Sigma_{\tau^{\beta}_{n, \delta}}\right) \frac{\log (3/ \delta)}{n}} + \frac{8}{3}\tau^{\beta}_{n, \delta}\sqrt{\overline{k}}\frac{\log (3/ \delta)}{n},
\end{equation*}
where $\Sigma_s$ is defined in Lemma~\ref{lemma:exact-expression-second-moment-covariances}.
\end{theorem}

\begin{example}[Countable state time-homogeneous Markov chains]
Irreducible, aperiodic and stationary countable Markov chains are always $\beta$-mixing \citep{bradley2005basic}.
Assume that there exist $\beta_1 \in \bbR_+$ and $b \in (1, \infty)$ such that for any $t\in \bbN, \beta(t) \leq  \beta_1/t^{b}$, that is, the chain is only known to be algebraically mixing.
    Interestingly, such a chain may not be $\phi$-mixing, thus this example cannot be recovered from Theorem~\ref{theorem:variance-aware-mmd-estimation-phi-mixing}.
\end{example}

\section{Applications}
\label{section:applications}

We put the apparatus developed  in Section~\ref{section:variance-aware} and Section~\ref{section:variance-empirical} to application in the context of hypothesis testing and robust parametric estimation. In this section, all considered kernels will be TI unless otherwise specified, and $X_1,\dots,X_n$ are i.i.d. from $\bbP \in  \calP(\calX)$. We  introduce a statistical model $\calM = \set{\bbP_\theta \colon \theta \in \Theta}$ indexed by the parameter space $\Theta$. Examples of models studied in the literature on MMD include parametric models such as the Gaussian model $\bbP_\theta = \calN(\theta, \sigma^2 I_d)$ (with $\sigma^2$ known) and mixture of Gaussians \citep{briol2019statistical,cherief2022finite}, copulas~\cite{alquier2022estimation} but also more complex models such as generative adversarial networks~\citep{dziugaite2015training,li2015generative}, stochastic volatility models and stochastic differential equations~\citep{briol2019statistical}.
Following \citet{briol2019statistical},
our estimator will be the closest element\footnote{
Note that it might be that the infimum is reached by multiple elements of $\Theta$, or is not reached. The first case does not lead to any difficulty, we must simply define a rule to break ties (for example, we can equip $\Theta$ with a total order and chose the smallest minimizer according to this order). While \citet{briol2019statistical} provide sufficient conditions to ensure that the infimum is reached, the non-existence of a minimizer is also not a problem in practice: all the non-asymptotic results in \citet{briol2019statistical} and \citet{cherief2022finite} can easily be extended to any $\epsilon$-minimizer, for $\epsilon$ small enough.
}
of the model $\calM$ to the empirical measure obtained from $X$, where the distance is measured in the RKHS,
\begin{equation}
\label{eq:parametric-estimator-definition}
    \nrm{\mu_{\bbP_{\widehat{\theta}_n(X_1,\dots, X_n)}} - \widehat{\mu}_{\bbP}(X_1, \dots, X_n)}_{\calH_k} = \inf_{\theta \in \Theta} \nrm{\mu_{\bbP_{\theta}} - \widehat{\mu}_{\bbP}(X_1, \dots, X_n)}_{\calH_k}.
\end{equation}
The computation of $\widehat{\theta}_n(X_1,\dots, X_n)$ is usually done via stochastic gradient descent and variants, see Remark~\ref{remark:observations-on-v} above.

\subsection{Hypothesis Testing}
\label{section:hypothesis-testing}

In this subsection we study hypothesis testing based on MMD. In the literature, two kinds of tests were proposed and studied: two-sample tests, and goodness of fit. In the two-sample test problem, we are given $X_1,\dots,X_n$ iid from some $\bbP_X$ and $Y_1,\dots,Y_m$ iid from some $\bbP_Y$, and we want to test $\mathbf{H}_0: \bbP_X = \bbP_Y$ against the alternative $\mathbf{H}_1: \bbP_X \neq \bbP_Y$. In the goodness of fit problem, we are given $X_1,\dots,X_n$ iid from $\bbP$ and we wish to test the hypothesis $\mathbf{H}_0: \bbP\in \calM =  \{\bbP_{\theta} \colon \theta\in\Theta \}$, against the alternative $\mathbf{H}_1: \bbP\notin\{\bbP_{\theta} \colon \theta\in\Theta \}$.

In the MMD literature, two-sample tests are more prevalent~\citep{gretton2009fast,gretton2012kernel}. Recent work  tackle goodness of fit testing \citep{pmlr-v48-chwialkowski16,jitkrittum2017linear}, albeit with an asymptotic treatment. We will study both problems here. In both cases, we propose a non-asymptotic treatment. That is, the level of the test is smaller than $\alpha$ for a finite sample size, and not only asymptotically. Moreover, in both cases, our data-dependent bound (Theorem~\ref{theorem:empirical-variance-mmd-bound}) allows to increase the power of the test, when compared to the procedure that would be based on the McDiarmid-based bound in~\eqref{eq:tolstikhin-bound}.

\subsubsection{Goodness-of-Fit Test}

Recall that we define the significance level $\alpha \in [0,1]$ of a test as the probability of outputting $\mathbf{H}_1$ when $\mathbf{H}_0$ is true.
Taking advantage of a non-asymptotic bound $\calB(X_1,\dots,X_n,\alpha)$ (for example  $\calB(X_1,\dots,X_n,\alpha)=\widehat{\calB}_{\alpha}$ given in Theorem~\ref{theorem:empirical-variance-mmd-bound}), we can design a test with prescribed level $\alpha$ for any $n$ as follows. We introduce the test statistic
$$ T(X_1, \dots, X_n) \eqdef \inf_{\theta\in\Theta} \nrm{\widehat{\mu}_{\bbP} - \mu_{\bbP_\theta}}_{\calH_k}, $$
and reject $\mathbf{H}_0$ on the ``critical set"
$$ \mathcal{C}(X_1, \dots, X_n) \eqdef \set{T > \calB(X_1,\dots,X_n,\alpha)}. $$
We show that the probability of rejection under the null hypothesis is at most $\alpha$ and that the test is consistent.

\begin{theorem}
\label{theorem:hypothesis-testing}
Let $n \in \N$, $X_1, \dots, X_n \sim \bbP$, and $k$ a translation invariant kernel. Let $\calB(\cdot,\alpha)$ be such that
\begin{enumerate}
\item[(a)] for any $\alpha$, with probability at least $1-\alpha$, $  \nrm{\ekme - \kme}_{\calH_{k}} \leq \calB(X_1,\dots,X_n,\alpha) $.
\end{enumerate}
Then
\begin{enumerate}
    \item[(i)] $\mathbb{P}_{\mathbf{H}_0} \left( \mathcal{C}(X_1, \dots, X_n) \right) \leq \alpha $.
\end{enumerate}
If, moreover,
\begin{enumerate}
\item[(b)] for any fixed $\alpha$, $\calB(X_1,\dots,X_n,\alpha) \xrightarrow[n\rightarrow\infty]{\text{a.s.}}  0$,
\item[(c)] there is a sequence $\alpha_n\rightarrow 0$ such that $\calB(X_1,\dots,X_n,\alpha_n) \xrightarrow[n\rightarrow\infty]{\text{a.s.}}  0$,
\end{enumerate}
then
\begin{enumerate}
  \item[(ii)] when the model $\calM$ is closed with respect to the MMD metric \footnote{For any $\mathbb{Q}$, if there is a sequence $(\theta_h)_{h\in\mathbb{N}}$ of elements of $\Theta$ such that $ \nrm{\mu_{\mathbb{Q}}- \mu_{\bbP_{\theta_h}}}_{\calH_k}  \xrightarrow[h\rightarrow\infty]{} 0  $, then $\mathbb{Q}=\bbP_{\theta}$ for some $\theta\in\Theta$.},
    $$\lim_{n \to \infty }\PR[\mathbf{H}_1]{\mathcal{C}(X_1, \dots, X_n)} = 1,$$ more precisely,
    $$ 1-\PR[\mathbf{H}_1]{\mathcal{C}(X_1, \dots, X_n)} = \mathcal{O}(\alpha_n). $$
\end{enumerate}
\end{theorem}
Obviously, the empirical bound $\widehat{\calB}_{k, \alpha}$ of Theorem~\ref{theorem:empirical-variance-mmd-bound} satisfies $(a)$ and $(b)$ (regardless of $\bbP$), and $(c)$ with $\alpha_n = \exp(-n^{1-\epsilon}) $ for any fixed $\epsilon>0$. So does the McDiarmid-based bound in~\eqref{eq:tolstikhin-bound}. On the other hand, we claimed that the empirical bound is smaller than bounds that does not take the variance into account. In other words, the power of the test $\PR[\mathbf{H}_1]{\mathcal{C}(X_1, \dots, X_n)} $ for a finite $n$ will be larger if we use the variance aware bound.

\begin{example}
We first consider a single hypothesis test (that is a special case of goodness-of-fit when $\Theta$ is a singleton). We consider data in $\mathbb{R}^2$: here $\bbP = \mathcal{N}(0,\sigma^2 I_2) $ with $\bbP_\theta = \mathcal{N}(\theta, I_2)$, and $\Theta = \{(1,1)\}$. In other words, $\mathbf{H}_0$ is true iff $\sigma=1$.  Note that in this case, if $Y_1, \dots, Y_n$ are iid from $\bbP_\theta=\bbP_{\{(1,1\}}$, then $q_{1-\alpha}$ defined as the $(1-\alpha)$-quantile of $\nrm{\mu_{\bbP_{0}} - \widehat{\mu}_{\bbP}(Y_1, \dots, Y_n)}_{\calH_k}$ allows to define a test with rejection zone given by $\mathcal{C}'(X_1,\dots,X_n) = \{ T > q_{1-\alpha} \} $ that satisfies $\mathbb{P}_{\mathbf{H}_0} \left( \mathcal{C}'(X_1, \dots, X_n) \right) = \alpha $ by definition. It is a very natural procedure to estimate $q_{1-\alpha}$ by Monte-Carlo, by sampling multiple times $Y_1, \dots, Y_n$ from $\bbP_\theta=\bbP_{\{(1,1\}}$. This leads to a Monte-Carlo estimator $\hat{q}_{1-\alpha}$ of $q_{1-\alpha}$.
In our simulations, we sample $X_1,\dots,X_n$ from $\bbP$, perform the test based on $\hat{q}_{1-\alpha}$, the test based on the empirical bound and the test based on McDiarmid bound . This is repeated $100$ times for each value $\sigma\in\{0,1/50,2/50,\dots,1\}$. We report the frequency of rejections in Figure~\ref{figure:test:simple}. The kernel used is a Gaussian kernel with $\gamma=1$, and we consider sample sized $n\in\{16,40,100,250\}$. We observe that, when compared to test based $\hat{q}_{1-\alpha}$, both tests based on bounds have a weak power (note that we did not try to optimize $\gamma$ for now, this question will be tackled later). However, the test based on the empirical bound indeed rejects more often $\mathbf{H}_0$ when $\mathbf{H}_1$ is true. The improvement is clearer for small sample sizes.
\begin{figure}
\begin{center}
    \includegraphics[width=0.49\textwidth]{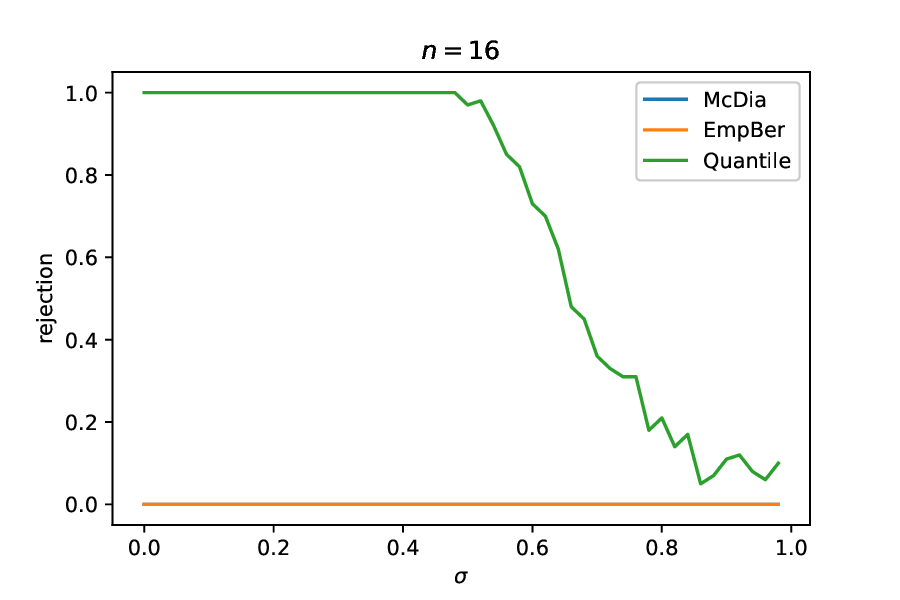}
    \includegraphics[width=0.49\textwidth]{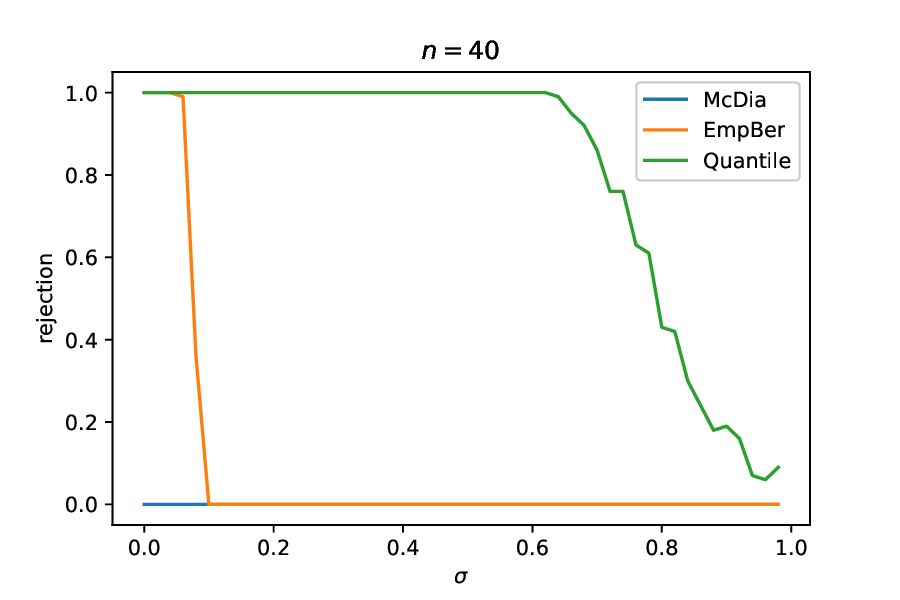}
    \includegraphics[width=0.49\textwidth]{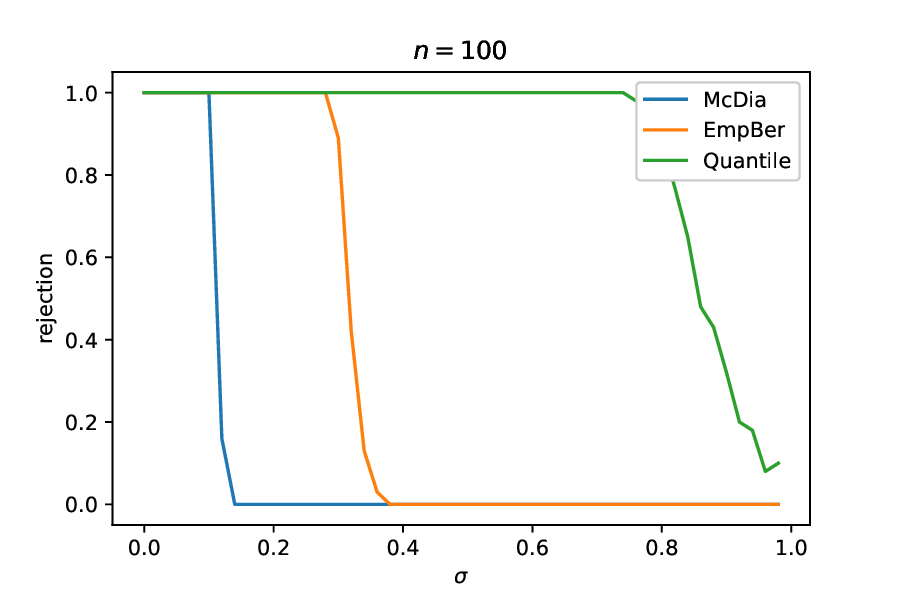}
    \includegraphics[width=0.49\textwidth]{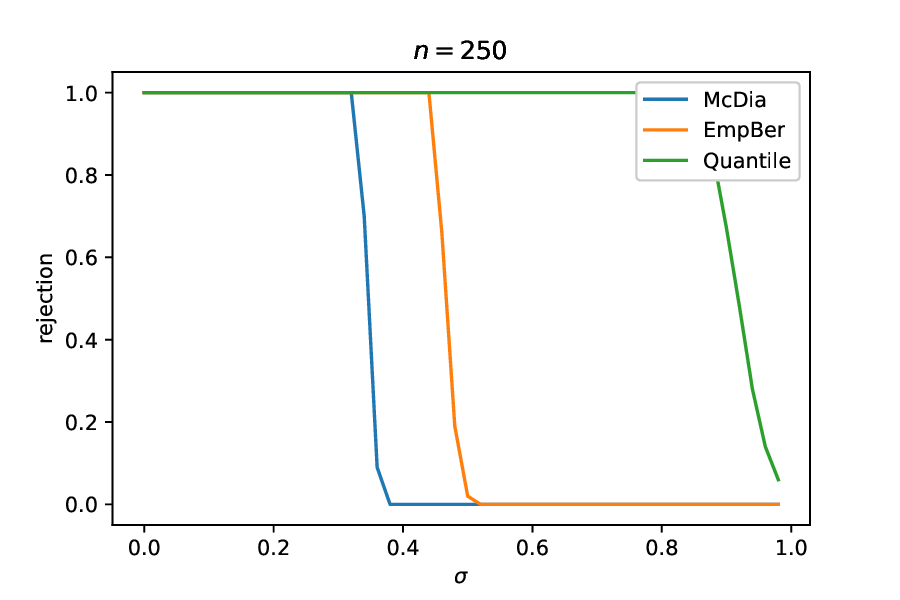}
    \caption{Comparison of the test based on the Bernstein empirical (EmpBer) bound, versus the test based on McDiarmid bound (McDia), and the test based on the Monte-Carlo estimation of the quantile $q_{1-\alpha}$. Frequency of rejection of $\mathbf{H}_0:\bbP\in\{\mathcal{N}((1,1), I_2)\}$ as a function of $\sigma$ with $\bbP= \mathcal{N}(0,\sigma^2 I_2) $. 
    }
    \label{figure:test:simple}
\end{center}
\end{figure}
\end{example}

\begin{example}
We now consider a proper goodness-of-fit test in $\mathbb{R}^2$: we still consider $\bbP = \mathcal{N}(0,\sigma^2 I_2) $ and $\bbP_\theta = \mathcal{N}(\theta, I_2)$, with $\Theta=\mathbb{R}^2$.
It is important to observe that in this case, because $\mathbf{H}_0$ is composite, we don't have a natural definition for $q_{1-\alpha}$ as in the previous example.
We sample $X_1,\dots,X_n$ from $\bbP$, perform the test based on the empirical bound and the test based on McDiarmid bound.
This is repeated $100$ times for each value $\sigma\in\{0,1/50,2/50,\dots,1\}$. We report the frequency of rejections in Figure~\ref{figure:test}. The kernel used is a Gaussian kernel with $\gamma=1$, and we consider sample sized $n\in\{16,40,100,250\}$. Similar comments to the previous case apply, but in this case, this makes the test based on the empirical Bernstein bound the best test available.
\begin{figure}
\begin{center}
    \includegraphics[width=0.49\textwidth]{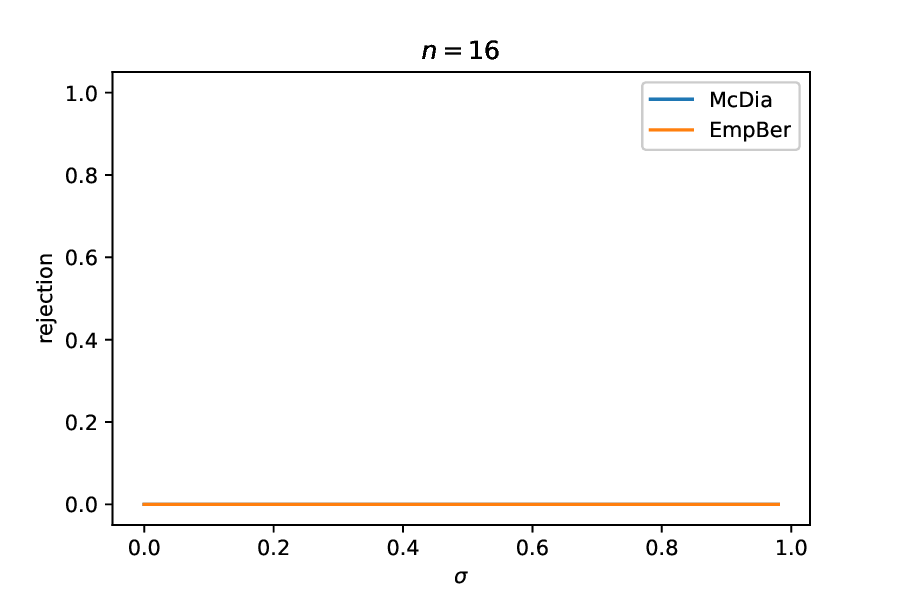}
    \includegraphics[width=0.49\textwidth]{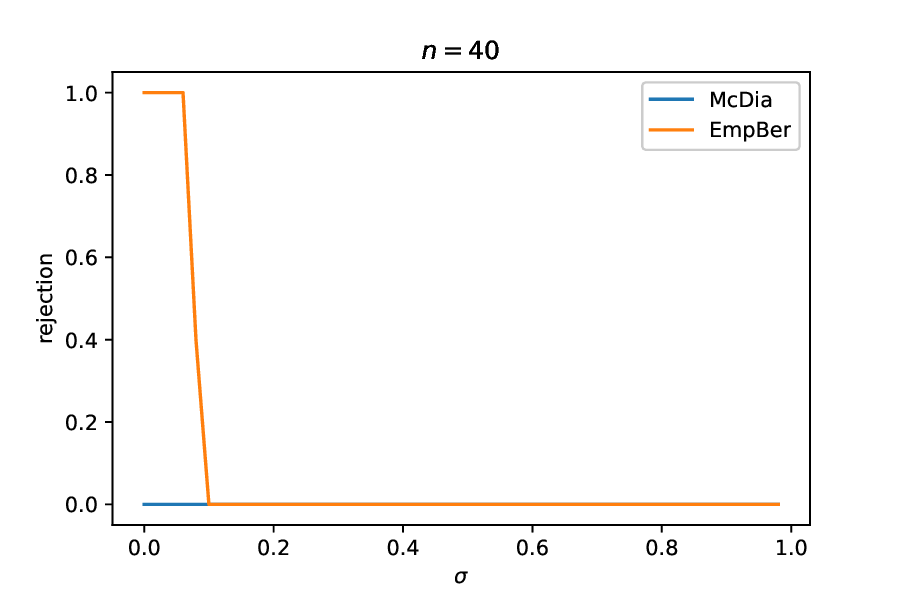}
    \includegraphics[width=0.49\textwidth]{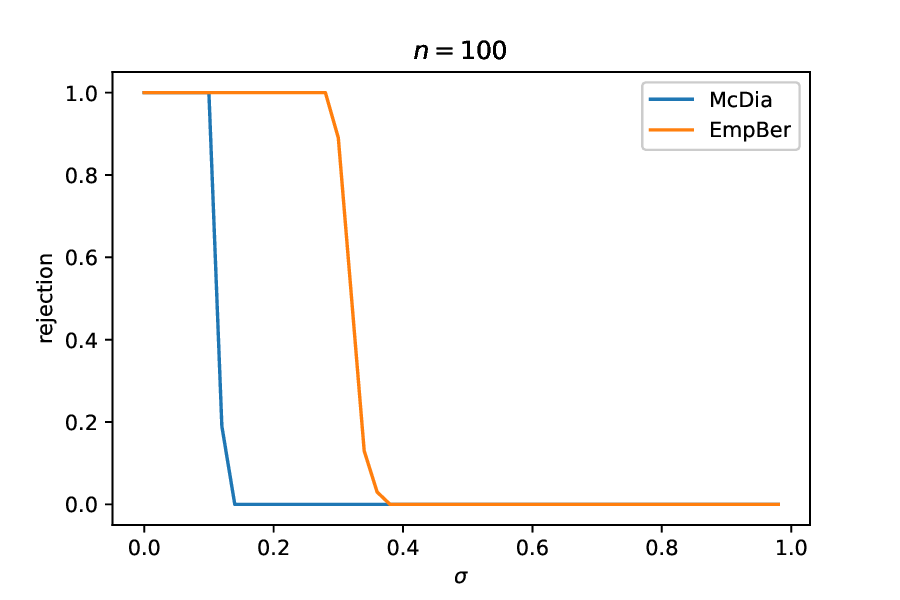}
    \includegraphics[width=0.49\textwidth]{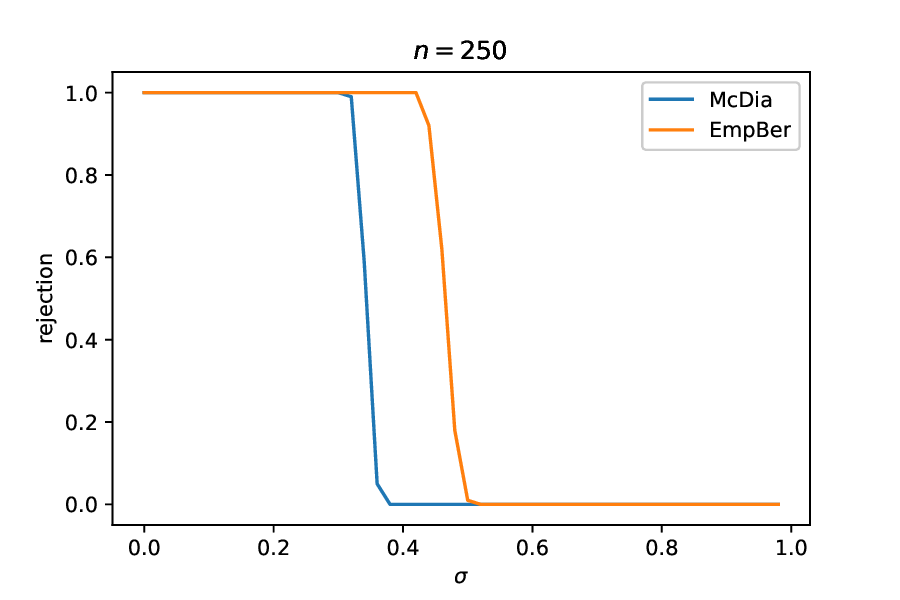}
    \caption{Comparison of the test based on the Bernstein empirical (EmpBer) bound, versus the test based on McDiarmid bound (McDia). Frequency of rejection of $\mathbf{H}_0:\bbP\in\{\mathcal{N}(\theta, I_2),\theta\in\mathbb{R}^2\}$ as a function of $\sigma$ with $\bbP= \mathcal{N}(0,\sigma^2 I_2) $. 
    }
    \label{figure:test}
\end{center}
\end{figure}
\end{example}
\subsubsection{Two-Sample Test}

Let $(n,m) \in \N^2$, $X_1, \dots, X_n \sim \bbP_X$ and $Y_1, \dots, Y_m \sim \bbP_Y$. Here, we wish to test $\mathbf{H}_0: \bbP_X = \bbP_Y$ against $\mathbf{H}_1: \bbP_X \neq \bbP_Y$.
This time, we use the test statistic
$$ T_2(X_1, \dots, X_n,Y_1,\dots,Y_m) \eqdef \nrm{\widehat{\mu}_{\bbP}(X_1,\dots,X_n) - \widehat{\mu}_{\bbP}(Y_1,\dots,Y_m)}_{\calH_k}, $$
and reject $\mathbf{H}_0$ on the ``critical set"
$$ \mathcal{C}_2(X_1, \dots, X_n,Y_1,\dots,Y_m) \eqdef \set{T_2 >
\calB(X_1,\dots,X_n,\alpha/2)+\calB(Y_1,\dots,Y_m,\alpha/2)}. $$

\begin{theorem}
\label{theorem:hypothesis-testing_2}
Let $(n,m) \in \N^2$, $X_1, \dots, X_n \sim \bbP_X$ and $Y_1, \dots,Y_m \sim \bbP_Y$, and $k$ a translation invariant kernel. Let $\calB(\cdot,\alpha)$ be such that
\begin{enumerate}
\item[(a)] for any $\alpha$, each of these inequalities hold with probability at least $1-\alpha$ each:
\begin{align*}  \nrm{\widehat{\mu}_{\bbP}(X_1,\dots,X_n) - \mu_{\bbP_X}  }_{\calH_k} & \leq \calB(X_1,\dots,X_n,\alpha) \\ \nrm{ \widehat{\mu}_{\bbP}(Y_1,\dots,Y_m)-\mu_{\bbP_Y} }_{\calH_k} & \leq \calB(Y_1,\dots,Y_m,\alpha).
\end{align*}
\end{enumerate}
Then
\begin{enumerate}
\item[(i)] $\mathbb{P}_{\mathbf{H}_0} \left( \mathcal{C}_2(X_1, \dots, X_n,Y_1,\dots,Y_m) \right) \leq \alpha $.
\end{enumerate}
Moreover, if
\begin{enumerate}
\item[(b)] for any fixed $\alpha$, $\calB(X_1,\dots,X_n,\alpha) \xrightarrow[n\rightarrow\infty]{\text{a.s.}}  0$ and $\calB(Y_1,\dots,Y_n,\alpha) \xrightarrow[m\rightarrow\infty]{\text{a.s.}}  0$,
\item[(c)] there is a sequence $\alpha_n\rightarrow 0$ such that $$\calB(X_1,\dots,X_n,\alpha_n) \xrightarrow[n\rightarrow\infty]{\text{a.s.}}  0 \text{ and } \calB(Y_1,\dots,Y_n,\alpha_m) \xrightarrow[m\rightarrow\infty]{\text{a.s.}}  0,$$
\end{enumerate}
then
\begin{enumerate}
    \item[(ii)] $\lim_{n,m \to \infty }\PR[\mathbf{H}_1]{\mathcal{C}_2(X_1, \dots, X_n,Y_1,\dots,T_m)} = 1.$
\end{enumerate}
\end{theorem}
Here again, the empirical bound $\widehat{\calB}_{k, \alpha}$ of Theorem~\ref{theorem:empirical-variance-mmd-bound} satisfies $(a)$ and $(b)$, and $(c)$ with $\alpha_n = \exp(-n^{1-\epsilon}) $ for any $\epsilon>0$.

\subsection{Robust Parametric Estimation under Huber Contamination}
\label{section:parametric-estimation}
\citet[Proof~of~Theorem~3.1]{cherief2022finite} show how to upper-bound~\eqref{eq:parametric-estimator-definition} by the estimation error of the empirical measure in the RKHS and the approximation error of the model,
\begin{equation}
\label{equation:cherief-1}
    \nrm{\mu_{\bbP_{\widehat{\theta}_n}} - \mu_{\bbP}}_{\calH_k} \leq \inf_{\theta \in \Theta} \nrm{\mu_{\bbP_\theta} - \mu_{\bbP}}_{\calH_k} + 2 \nrm{\widehat{\mu}_{\bbP} - \mu_{\bbP}}_{\calH_k}.
\end{equation}
Recall that in the Huber contamination model \citep{huber2011robust}, the observations $X_1, \dots, X_n$ are drawn independently from the mixture
\begin{equation*}
    \bbP = (1 - \xi) \bbP_{\theta_0} + \xi \bbH,
\end{equation*}
where $\xi \in [0, 1/2)$ is the contamination rate, $\bbH$ is some unknown noise distribution and $\theta_0 \in \Theta$. In this setting, we can control the approximation error of the model: starting from \eqref{equation:cherief-1}, and adapting the argument of \citet[Corollary~3.3]{cherief2022finite}, we obtain
\begin{equation*}
    \nrm{\mu_{\bbP_{\widehat{\theta}_n}} - \mu_{\bbP_{\theta_0}}}_{\calH_k} \leq 4 \xi \sqrt{\overline{\psi}} + 2 \nrm{\widehat{\mu}_{\bbP} - \mu_{\bbP}}_{\calH_k}.
\end{equation*}
We wish to apply Theorem~\ref{theorem:variance-aware-mmd-estimation} to the second term.
However, when the sample is contaminated, the variance term $v(\bbP)$ we collect also depends on the properties of $\bbH$, which are unknown---for all we know $\bbH$ does not even have finite moments. Fortunately, we now see how to  bound $v(\bbP)$ is terms of $v(\bbP_{\theta_0})$ and $\xi$.

\subsubsection{Improved Confidence Bounds in the Huber Contamination Model}
\label{section:improving-bounds-contaminated-data}

\begin{lemma}
\label{lemma:variance-huber-contamination}
In the Huber contamination model, that is, $
    \bbP = (1 - \xi) \bbP_{\theta_0} + \xi \bbH$, writing $v = v(\bbP)$ and $v_0 = v(\bbP_{\theta_0})$, it holds that
\begin{equation}
\label{eq:huber-contamination-variance-bound-precise}
\begin{split}
    v \leq v_0 + 2 \xi(\Delta \psi - v_0) + \xi^2(v_0 - \Delta \psi).
\end{split}
\end{equation}    
The bound can be streamlined; when $\psi \geq 0$,
\begin{equation*}
    v \leq (1 - 2 \xi) v_0 + 2 \xi \overline{\psi},
\end{equation*}
and otherwise
\begin{equation*}
    v \leq (1 - 2 \xi) v_0 + 2 \xi \Delta \psi.
\end{equation*}
\end{lemma}

\begin{remark}
In \eqref{eq:huber-contamination-variance-bound-precise},
observe that when $\psi \geq 0$ and  $\xi \to 1$, the first term vanishes, and $v \to \overline{\psi}$, recovering the distribution independent rates. 
Conversely, when $\xi \to 0$, we confirm that $v \to v_0$, which could be further bounded using properties of the model. In practice, we will focus on small values $\xi \ll 1/2$, and the simplified bounds will be sufficient.
\end{remark}
An application of Theorem~\ref{theorem:variance-aware-mmd-estimation} yields the following.
\begin{corollary}
\label{cor:huber}
In the Huber contamination model $
    \bbP = (1 - \xi) \bbP_{\theta_0} + \xi \bbH$, writing $v = v(\bbP)$ and $v_0 = v(\bbP_{\theta_0})$, with probability at least $1-\delta$,
\begin{equation*}
    \nrm{\mu_{\bbP_{\widehat{\theta}_n}} - \mu_{\bbP_{\theta_0}}}_{\calH_k} \leq 4 \xi \sqrt{\overline{\psi}} + 2 \sqrt{2 [(1 - 2 \xi) v_0 + 2 \xi \overline{\psi}] \frac{\log (2/ \delta)}{n}} + 2\frac{4 \sqrt{\overline{\psi}}}{3}\frac{\log (2/ \delta)}{n}.
\end{equation*}
\end{corollary}
When $\overline{\psi} = 1, v_0 \ll 1, \xi \ll 1$, the bound offers a significant improvement over \citet[Corollary~3.4]{cherief2022finite}.

\subsubsection{Improved Confidence Bounds in the Parameter Space}

In this subsection we still employ the Gaussian kernel $k_\gamma$ \eqref{definition:gaussian-kernel}.
It is instructive to analyze how robust estimation with respect to the MMD distance translates into what happens in the space of parameters. In fact, since we have freedom over the lengthscale parameter, we may understandably want to select $\gamma$ such that the distance in the space of parameters is kept small as well.

\begin{definition}[Link function $F$]
Let $\calM = \set{\bbP_\theta \colon \theta \in \Theta}$ indexed by the parameter space $\Theta\subset \calP(\R^d)$.
We say that a non-decreasing function $F_k \colon \R_+ \to \R_+$ is a link function for the model $\calM$ and kernel $k$ when for any $\theta, \theta' \in \Theta$,
\begin{equation*}
    \nrm{\theta - \theta'}_2 \leq F_k\left( \nrm{ \mu_{\bbP_\theta} - \mu_{\bbP_{\theta'}} }_{\calH_{k}} \right).
\end{equation*}
\end{definition}
\begin{remark}
Note that a link function $F_k$ will only provide nontrivial information if $F_k(0)=0$ and $F_k$ is continuous at $0$. The existence of such a nontrivial link function implies that the model is identifiable: if $\mu_{\bbP_\theta} = \mu_{\bbP_{\theta'}} $ then $0=F_k(\nrm{ \mu_{\bbP_\theta} - \mu_{\bbP_{\theta'}} }_{\calH_{k}}) \geq \nrm{\theta-\theta'}_{2} $, hence $\theta=\theta'$. In this case, there is a unique function $p:\mathcal{M}\rightarrow\Theta $ such that $p(\bbP_\theta)=\theta$, and $F_k$ is a modulus of continuity of $p$.
\end{remark}
A direct application of Corollary~\ref{cor:huber} gives the following.
\begin{corollary}
\label{cor:huber-parameter}
In the setting of Corollary~\ref{cor:huber}, assume that $\mathcal{M}$ and $k$ are such that a link function $F_k$ exists. With probability at least $1-\delta$,
\begin{equation*}
    \nrm{\widehat{\theta}_n - \theta_0 }_{2} \leq F_k \left( 4 \xi \sqrt{\overline{\psi}} + 2 \sqrt{2 [(1 - 2 \xi) v_0 + 2 \xi \overline{\psi}] \frac{\log (2/ \delta)}{n}} + 2\frac{4 \sqrt{\overline{\psi}}}{3}\frac{\log (2/ \delta)}{n} \right).
\end{equation*}
\end{corollary}

\emph{Gaussian location model, continued.}
We continue here Example~\ref{example:gaussian-known-variance}.
A direct computation (see for instance \citealt{cherief2022finite}) shows that
\begin{equation*}
    \nrm{\mu_{\bbP_{\theta}} - \mu_{\bbP_{\theta'}}}_{\calH_{k_\gamma}}^2 = 2 \left( \frac{\gamma^2}{2 \sigma^2 + \gamma^2}\right)^{d/2} \left[ 1 - \exp \left(- \frac{\nrm{\theta' - \theta}^2_2}{4 \sigma^2 + 2 \gamma^2} \right) \right].
\end{equation*}
In other words, a link function for the Gaussian location model is explicitly given by
\begin{equation*}
     F^2_{k_\gamma}(h) = -2(2 \sigma^2 + \gamma^2)  \log \left(1 - \frac{h}{2}\left( 1 + 2 \frac{ \sigma^2 }{\gamma^2}\right)^{d/2}  \right).
\end{equation*}
An application of Corollary~\ref{cor:huber-parameter}, and setting $\gamma = \lambda \sigma \sqrt{d}$, for some $\lambda >0$, proves that, with probability at least $1-\delta$,
\begin{equation}
\label{eq:bound-theta-huber-contamination}
\begin{split}
    \nrm{\theta - \widehat{\theta}_n}^2_2 &\leq -2\sigma^2 (2  + d \lambda^2) \log \Bigg\{ 1 -  \left(1 + \frac{2 }{d \lambda^2} \right)^{d/2} \Bigg(4 \xi  \\ & \qquad \qquad +  \sqrt{2 \left[(1 - 2 \xi)\left(1 - e^{-\frac{1}{\lambda^2}}\right) + 2 \xi \right]\frac{\log (2/ \delta)}{n}} + \frac{4}{3} \frac{\log (2/ \delta)}{n} \Bigg)^2 \Bigg\}
    \\
    & \eqdef G_{d,\sigma,n}(\lambda,\xi).
\end{split}
\end{equation}
For comparison, we provide the bound of~\citet[Proposition~4.1-Equation~2]{cherief2022finite} that is obtained by plugging~\eqref{eq:tolstikhin-bound} into the link function:
\begin{equation}
\begin{split}
\label{eq:bound-cherief}
\nrm{\theta - \widehat{\theta}_n}^2_2 &\leq -2\sigma^2 (2  + d \lambda^2) \log \set{ 1 -  \left(1 + \frac{2 }{d \lambda^2} \right)^{d/2} \left(4 \xi +  \frac{1+\sqrt{2\log (1/\delta)}}{\sqrt{n}}  \right)^2 } \\
&\eqdef H_{d,\sigma,n}(\lambda,\xi) .
\end{split}
\end{equation}
\begin{figure}
\begin{center}
    \includegraphics[width=0.49\textwidth]{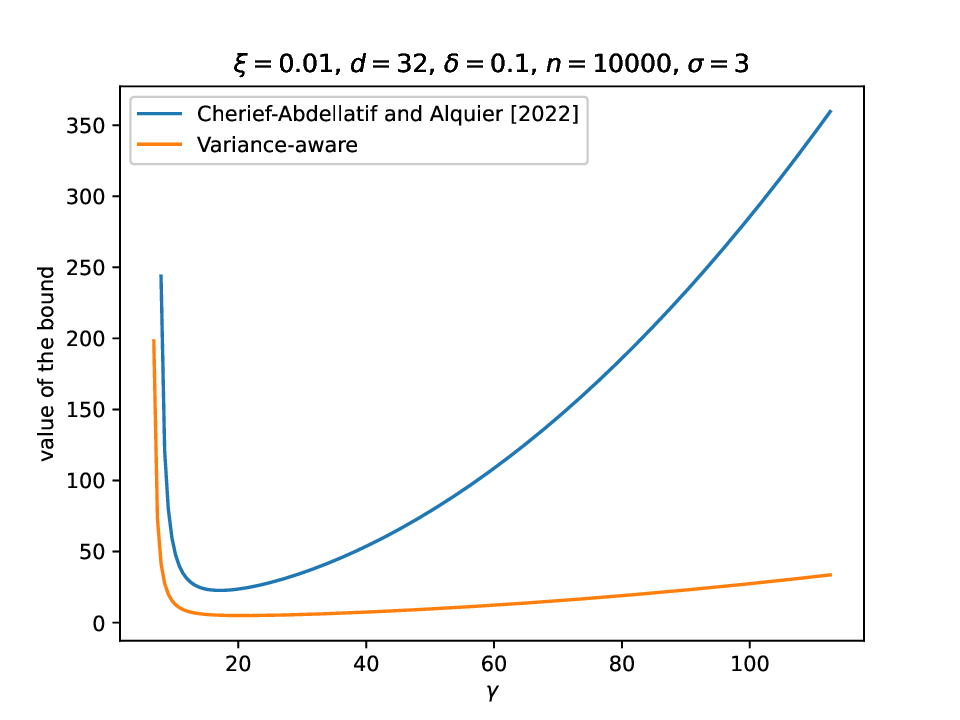}
    \includegraphics[width=0.49\textwidth]{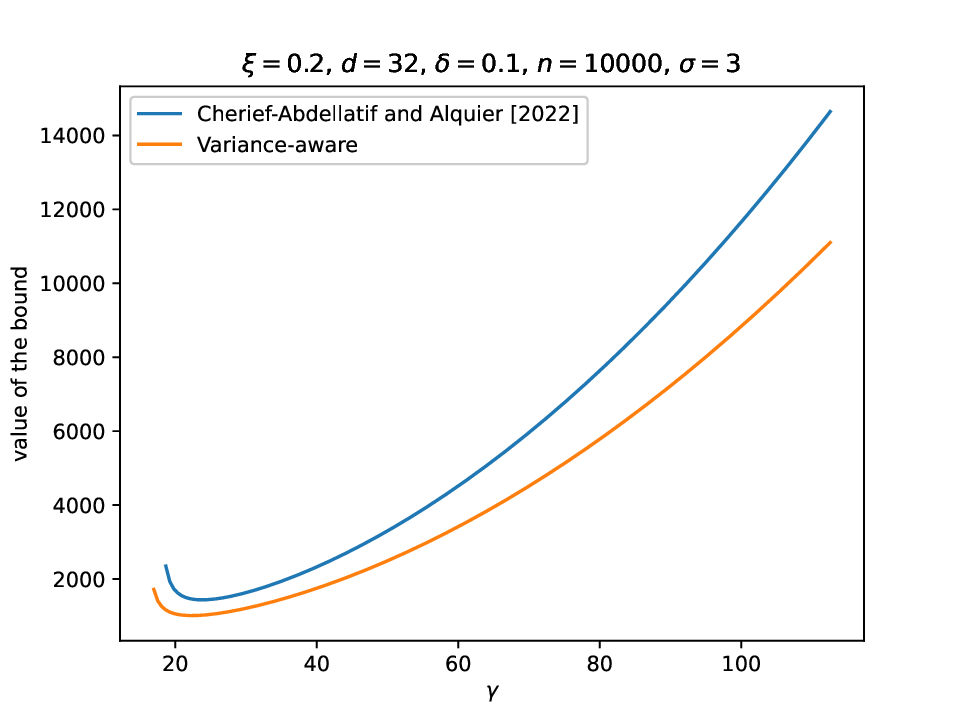}
    \caption{Comparison of the bounds in \eqref{eq:bound-theta-huber-contamination} versus \citet[Proposition~4.1-Equation~2]{cherief2022finite} as a function of $\gamma$.}
    \label{figure:experiment-alpha}
\end{center}
\end{figure}
\begin{figure}
\begin{center}
    \includegraphics[width=0.7\textwidth]{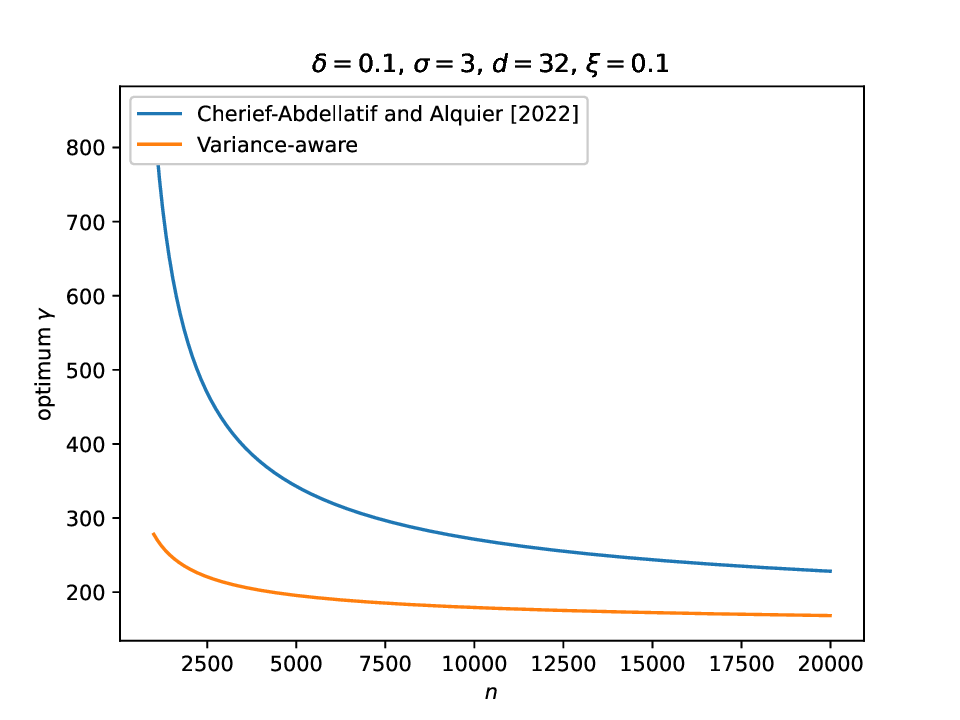}
    \caption{Comparison of \citet{cherief2022finite} and \eqref{eq:bound-theta-huber-contamination} for the optimal hyper-parameter $\gamma$ as a function of $n$.}
    \label{figure:experiment-small-sample}
\end{center}
\end{figure}
In order to compare the tightness of the bounds in Eq. \eqref{eq:bound-cherief} and Eq. \eqref{eq:bound-theta-huber-contamination}, 
we run experiments for fixed sample size ($n = 10000$), fixed confidence level ($\delta = 0.1$), fixed variance parameter ($\sigma = 3$) and two different contamination levels ($\xi = 0.01$ and $\xi = 0.2$). 
We compute and plot the bound obtained when varying the scale parameter $\gamma$ (see Figure~\ref{figure:experiment-alpha}).
For a contamination level $\xi = 0.1$, we also compare the quantitative decay of the two bounds---optimized for $\gamma$---as the sample size $n$ grows and report our results in \ref{figure:experiment-small-sample}.
Based on our experiments,
we make the following observations.
\begin{enumerate}[$(i)$]
    \item The new bound is always tighter than the one of \citet{cherief2022finite}.
    \item When the contamination level increases, the two bounds are getting closer together. This is expected since we are losing the benefit of the variance factor.
    \item Especially under weak contamination, the new bound is much flatter in the sense where overshooting for the value of $\gamma$ does not lead to a catastrophic degradation of the bound.
    \item The new bound performs exceptionally well in the small sample setting.
\end{enumerate}

We conclude this section by showing heuristically that the bound in~\eqref{eq:bound-theta-huber-contamination} can always be made smaller than the bound in~\eqref{eq:bound-cherief}, at least asymptotically in $n$, with an adequate choice of $\lambda$. For the sake of simplicity, we work in the non-contaminated setting $\xi=0$.
First,
\begin{equation*}
H_{d,\sigma,n}(\lambda,0)  = 4 \sigma^2 (2  + d \lambda^2)\left(1 + \frac{2 }{d \lambda^2} \right)^{d/2} \frac{ \left(1+\sqrt{2\log (1/ \delta)} \right)^2}{ n} (1+o(1)),
\end{equation*}
and the first order term is exactly minimized for $\lambda=1$, it leads to
\begin{equation*}
H_{d,\sigma,n}(1,0)  = 4 \sigma^2 (2  + d) \left(1 + \frac{2 }{d} \right)^{d/2} \frac{ \left(1+\sqrt{2\log (1/ \delta)} \right)^2}{ n} (1+o(1)).
\end{equation*}

An exact minimization of the  bound in~\eqref{eq:bound-theta-huber-contamination} is not feasible, but we will propose a choice of $\lambda=\lambda_n$ that will lead to an improvement on $H_{d,\sigma,n}(1,0)$. First, let us assume that $\lambda_1>0$ and that $(\lambda_n)_{n \in \N}$ is a non-decreasing sequence. We obtain:
\begin{align*}
G_{d,\sigma,n}(\lambda_n,0) & = 4 \sigma^2 (2  + d \lambda_n^2)\left(1 + \frac{2 }{d \lambda_n^2} \right)^{d/2} \left( 1-e^{-\frac{1}{\lambda_n^2}} \right) \frac{ \log (2/ \delta)}{ n} (1+o(1))
\\
& = 4 \sigma^2 (2  + d \lambda_n^2)\left(1 + \frac{2 }{d \lambda_n^2} \right)^{d/2} \frac{ \log (2/ \delta)}{ \lambda_n^2 n} (1+o(1))
\\
& = 4 d \sigma^2 \left(1 + \frac{2 }{d \lambda_n^2} \right)^{1+d/2} \frac{\log (2/ \delta)}{  n} (1+o(1)).
\end{align*}
It is clear that choosing $\lambda_n\rightarrow\infty$ when $n\rightarrow\infty$ gives:
\begin{equation*}
G_{d,\sigma,n}(\lambda_n,0)  = 4 d \sigma^2 \frac{\log (2/ \delta)}{n} (1+o(1)).
\end{equation*}
Finally, observe that
\begin{equation*}
\frac{
H_{d,\sigma,n}(1,0) 
}
{
G_{d,\sigma,n}(\lambda_n,0)
}
= 2 \underbrace{\left(1 + \frac{2 }{d} \right)^{d/2+1}}_{ \xrightarrow[d\rightarrow\infty]{\searrow} e } \underbrace{ \frac{ (1+\sqrt{2\log 1/ \delta} )^2 } { (\sqrt{2\log (2/ \delta)} )^2 } }_{\geq 1 \text{ for } \delta \in (0, 9/10)} (1+o(1)) \geq 2 e (1+o(1)).
\end{equation*}
Unfortunately, this first order analysis is too crude to provide an accurate recommendation on the choice of $\lambda_n$. On the other hand, it shows that our variance-aware bound leads to a significant improvement over the bound of~\cite{cherief2022finite} for Gaussian mean estimation. Moreover, the choice $\lambda_n\rightarrow\infty$ is in accordance with the considerations on the asymptotic variance in~\cite{briol2019statistical} in this model.

\begin{figure}
\begin{center}
    \includegraphics[width=0.49\textwidth]{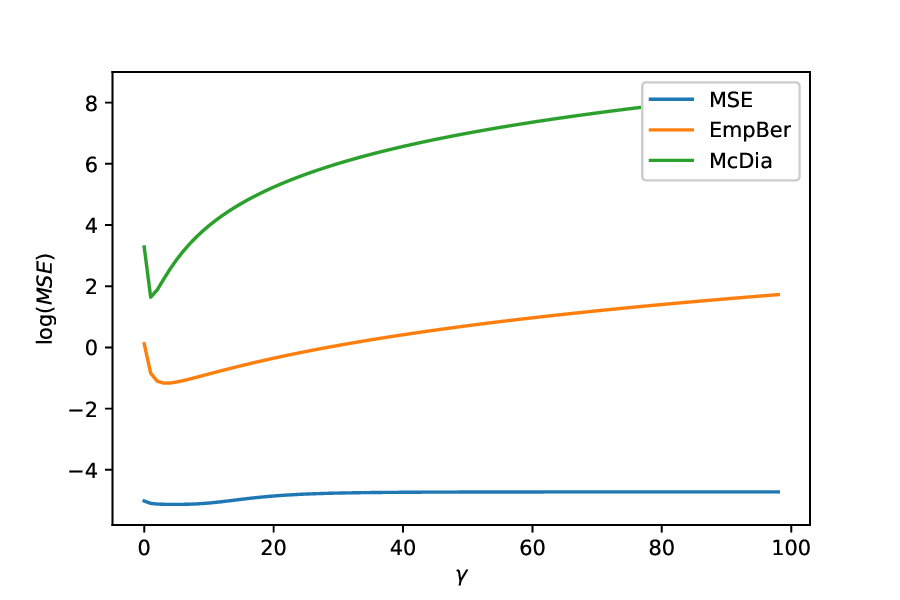}
    \includegraphics[width=0.49\textwidth]{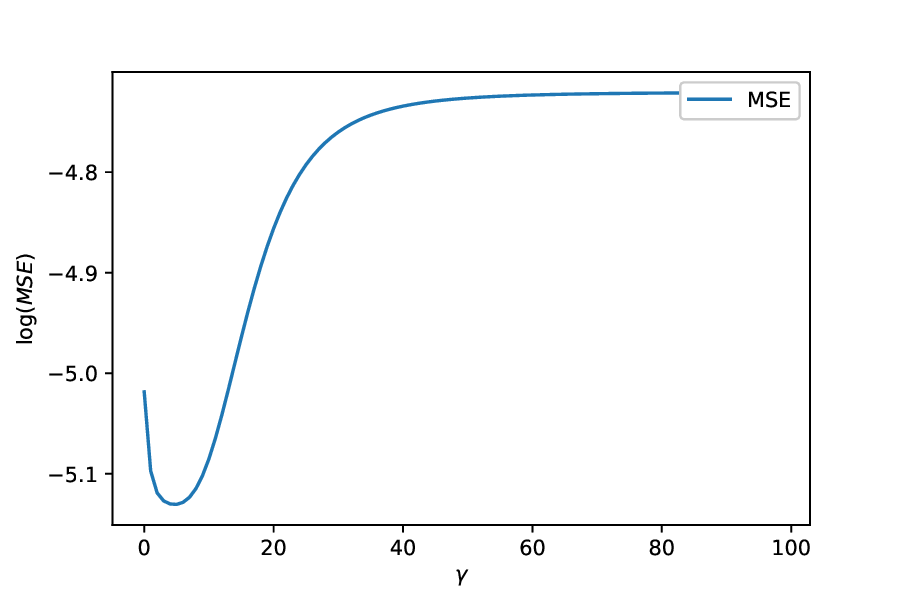}
    \caption{Comparison of the bounds in \eqref{eq:bound-theta-huber-contamination} versus \citet[Proposition~4.1-Equation~2]{cherief2022finite} as a function of $\gamma$, versus the empirical MSE on $50$ experiments, $n=500$, log-scale, $\xi=0$, $\delta=0.05$. Right panel: we zoom on the MSE.}
    \label{figure:mse}
\end{center}
\end{figure}

We illustrate this with an experiment in Figure~\ref{figure:mse}. We observe that: indeed, the variance-aware bound decrease with $\gamma\propto\lambda$, while the McDiarmid bound has a minimum. Interestingly, the true MSE also has a minimum. Note however how the MSE is actually very flat as a function of $\gamma$: in other words, in this experiment, the choice of $\gamma$ does not matter so much on the performance of the MMD estimator. The possibility to use the variance-aware bound to calibrate $\gamma$ in practice should be investigated further in other models.

\section{Proofs}
\label{section:proofs}
In this section, we provide the detailed technical proofs that were deferred in earlier parts of the paper.

\subsection{Proof of Theorem~\ref{theorem:variance-aware-mmd-estimation}}
In \citet{tolstikhin2017minimax}, the authors mention but do not pursue the idea that Bernstein's inequality in separable Hilbert
spaces yields the proper rate of $\bigO_P(n^{-1/2})$. We follow-up on their idea in this proof.
First recall that
\begin{equation*}
\begin{split}
    \ekme \eqdef \frac{1}{n}\sum_{t = 1}^{n} k(X_t, \cdot), \qquad \kme \eqdef \E[X \sim \bbP]{k(X, \cdot)},
\end{split}
\end{equation*}
and that we can write
\begin{equation*}
\begin{split}
    \nrm{\ekme - \kme}_{\calH_k}
    &=  \frac{1}{n} \nrm{\sum_{t = 1}^{n} Z_t}_{\calH_k},
\end{split}
\end{equation*}
where for $1 \leq t \leq n$,
\begin{equation*}
\begin{split}
    Z_t \eqdef k(X_t, \cdot) - \E[X \sim \bbP]{k(X, \cdot)},\\
\end{split}
\end{equation*}
is centered. For $\eps >0$, we are thus interested in bounding the probability of the following deviation,
\begin{equation*}
    \PR{ \max_{1 \leq s \leq n} \nrm{  \sum_{t = 1}^{s} Z_t}_{\calH_k} > n \eps }.
\end{equation*}
We rely on a Bernstein-type inequality in Hilbert spaces, credited to \citet{pinelis1985remarks} (see also \citealt[Theorem~3.3.4]{yurinsky1995sums}), which is reported in Theorem~\ref{theorem:bernstein-inequality-hilbert-space}.
We first observe that the RKHS $\calH_k$ is separable by separability of the topological space
$\calX$ and continuity of $k$.
In order to invoke Theorem~\ref{theorem:bernstein-inequality-hilbert-space}, we need to control from above the $\sum_{t = 1}^{n} \bbE \nrm{Z_t}_{\calH_k}^p$,
for any $p \geq 3$.
For $p=2$,
\begin{equation*}
    \sum_{t = 1}^{n} \bbE \nrm{k(X_t, \cdot) - \kme}_{\calH_k}^2 = n v.
\end{equation*}
Moving on to higher moments, for $p \geq 3$, we have that
\begin{equation*}
\begin{split}
    \nrm{Z_t}_{\calH_k} &\leq \nrm{k(X_t, \cdot)}_{\calH_k} + \nrm{\E[X \sim \bbP]{k(X, \cdot)}}_{\calH_k} \\
    &\leq \sqrt{k(X_t, X_t)} + \sup_{x \in \calX} \sqrt{k(x,x)} \\
    &\leq 2 \sqrt{\overline{k}}.
\end{split}
\end{equation*}
As a result,
\begin{equation*}
\begin{split}
    \sum_{t=1}^{n} \bbE \nrm{Z_t}_{\calH_k}^p &\leq   \left(\sum_{t=1}^{n} \bbE \nrm{Z_t}_{\calH_k}^2 \right) \left(2 
    \sqrt{\overline{k}}\right)^{p-2} \\ &\leq n v (2 
    \sqrt{\overline{k}})^{p-2} \leq \frac{p!}{2} n v \left(\frac{2 \sqrt{\overline{k}}}{3}\right)^{p-2}.
\end{split}
\end{equation*}
Invoking Theorem~\ref{theorem:bernstein-inequality-hilbert-space} with $G^2 = n v$, $H = 2 \sqrt{\overline{k}}/3$, we obtain that for any $\eps > 0$,
\begin{equation*}
\begin{split}
    \PR{ \max_{1 \leq s \leq n} \nrm{\sum_{t = 1}^{s} Z_t}_{\calH_k} > n \eps } &\leq 2\exp \left(- \frac{(\eps n)^2}{2(n v + (\eps n) 2  \sqrt{\overline{k}}/3)} \right) \\ &= 2 \exp \left(- \frac{n\eps^2}{2\left(v + 2\eps  \sqrt{\overline{k}}/3\right)} \right).
\end{split}
\end{equation*}
\qedsymbol

\subsection{Proof of Lemma~\ref{lemma:translation-invariant-weakly-self-bounding}}
We first verify the bounded-differences property.
\begin{equation*}
\begin{split}
 &\widehat{v}(x) - \inf_{x_{t_0}' \in \calX} \widehat{v}\left(x^{(t_0)}\right) \\
 =&  \overline{\psi} - \frac{1}{(n-1)n} \sum_{t \neq s} \psi (x_t - x_s) - \inf_{x_{t_0}' \in \calX} \left[\overline{\psi} - \frac{1}{(n-1)n} \sum_{t \neq s} \psi \left(x^{(t_0)}_t - x^{(t_0)}_s\right) \right]\\
  =& \frac{1}{(n-1)n} \left( \sup_{x_{t_0}'  \in \calX}  \set{ \sum_{t \neq s} \psi \left(x^{(t_0)}_t - x^{(t_0)}_s\right)} -  \sum_{t \neq s} \psi (x_t - x_s) \right) \\
  =& \frac{2}{(n-1)n} \left(  \sup_{x_{t_0}'  \in \calX}  \set{ \sum_{\substack{t \in [n] \\ t \neq t_0 }} \psi \left(x'_{t_0} - x_t\right)} -  \sum_{\substack{t \in [n] \\ t \neq t_0 }} \psi (x_{t_0} - x_t) \right) \\
  \leq& \frac{2}{(n-1)n} \left(   (n-1) \overline{\psi} -  \sum_{\substack{t \in [n] \\ t \neq t_0 }} \psi (x_{t_0} - x_t) \right) \\
  \stackrel{(\star)}{=}& \frac{2}{n} \left( \overline{\psi} - \frac{1}{n-1} \sum_{\substack{t \in [n] \\ t \neq t_0 }} \psi (x_{t_0} - x_t) \right) \\
  \leq& \frac{2}{n} \Delta \psi.
\end{split}
\end{equation*}
It follows from $(\star)$ that,
\begin{equation*}
\begin{split}
    \sum_{t_0 = 1}^{n} \left( \widehat{v}(x) - \inf_{x_{t_0}'  \in \calX} \widehat{v}\left(x^{(t_0)}\right) \right)^2 &\leq  \sum_{t_0 = 1}^{n} \left(\frac{2}{n} \left( \overline{\psi} - \frac{1}{n-1} \sum_{\substack{t \in [n] \\ t \neq t_0 }} \psi (x_{t_0} - x_t) \right)\right)^2 \\
    &= \frac{4}{n^2} \sum_{t_0 = 1}^{n} \left( \frac{1}{n-1} \sum_{\substack{t \in [n] \\ t \neq t_0 }} \left( \overline{\psi} -  \psi (x_{t_0} - x_t) \right)\right)^2 \\
    &\leq \frac{4}{n^2} \sum_{t_0 = 1}^{n} \frac{1}{n-1} \sum_{\substack{t \in [n] \\ t \neq t_0 }} \left(   \overline{\psi} -  \psi (x_{t_0} - x_t) \right)^2 \\
    &\leq \frac{4 \Delta \psi}{n} \sum_{t_0 = 1}^{n} \frac{1}{(n-1)n} \sum_{\substack{t \in [n] \\ t \neq t_0 }} \left(   \overline{\psi} -  \psi (x_{t_0} - x_t) \right) \\
    &= \frac{4 \Delta \psi}{n} \widehat{v}
\end{split}
\end{equation*}
where the second inequality is Jensen's.

\qedsymbol

\subsection{Proof of Lemma~\ref{lemma:concentration-square-root-v}}
From Lemma~\ref{lemma:translation-invariant-weakly-self-bounding}, the function
\begin{equation*}
    \calX^n \to \R, x \mapsto \frac{n}{2 \Delta \psi} \widehat{v}(x),
\end{equation*}
is weakly $(2, 0)$-self-bounding and has the bounded differences property in the sense of Definition~\ref{definition:weakly-self-bounding-and-bounded-differences-function}.
From the second statement of \citet[Theorem~13]{maurer2006concentration},
\begin{equation*}
    \PR{\E{\frac{n}{2 \Delta \psi} \widehat{v}(X)} - \frac{n}{2 \Delta \psi} \widehat{v}(X) > t} \leq \exp \left( - \frac{t^2}{4 \E{\frac{n}{2 \Delta \psi} \widehat{v}(X)}} \right),
\end{equation*}
that can be rewritten,
\begin{equation*}
    \PR{v -  \widehat{v}(X) > t} \leq \exp \left( - \frac{n t^2}{8 \Delta \psi v} \right).
\end{equation*}
In other words, with probability at least $1 - \delta$,
\begin{equation*}
    v - 2 \sqrt{v} \sqrt{ \frac{2 \Delta \psi \log (1/\delta)}{n}} \leq \widehat{v}(X).
\end{equation*}
Completing the square, and by sub-additivity of the square root function, we successively have with probability $1- \delta$,
\begin{equation*}
    \left(\sqrt{v} - \sqrt{\frac{2 \Delta \psi \log (1/\delta)}{n}}\right)^2 \leq \widehat{v}(X) + \frac{2 \Delta \psi \log (1/\delta)}{n}
\end{equation*}
\begin{equation*}
    \sqrt{v} - \sqrt{\frac{2 \Delta \psi \log (1/\delta)}{n}} \leq \sqrt{\widehat{v}(X) + \frac{2 \Delta \psi \log (1/\delta)}{n}}
\end{equation*}
\begin{equation*}
    \sqrt{v} \leq \sqrt{\widehat{v}(X)} + 2 \sqrt{\frac{2 \Delta \psi \log (1/\delta)}{n}},
\end{equation*}
which proves the lemma for $b = -1$.
On the other hand, from the first inequality in \citet[Theorem~13]{maurer2006concentration},
\begin{equation*}
    \PR{\frac{n}{2 \Delta \psi} \widehat{v}(X)  - \E{\frac{n}{2 \Delta \psi} \widehat{v}(X)} > t} \leq \exp \left( - \frac{t^2}{4 \E{\frac{n}{2 \Delta \psi} \widehat{v}(X)} + 2t} \right),
\end{equation*}
that can be rewritten as
\begin{equation*}
    \PR{\widehat{v}(X)  - v > t } \leq \exp \left( - \frac{n t^2}{8 \Delta \psi(v + t /2)} \right).
\end{equation*}
The positive solution $t_+$ for 
\begin{equation*}
    n t^2 = 8 \Delta \psi(v + t/2) \log (1/\delta)
\end{equation*}
is readily given by
\begin{equation*}
    t_+ =  \frac{2 \Delta \psi \log (1/ \delta)}{n} + \frac{4}{n} \sqrt{\left(\frac{\Delta \psi}{2} \log (1 /\delta)\right)\left(\frac{\Delta \psi}{2} \log (1 /\delta) + n v\right)}.
\end{equation*}
Thus, it holds with probability at least $1 - \delta$ that
\begin{equation*}
    \widehat{v}(X) \leq v + t_+,
\end{equation*}
and from sub-additivity of the square root, successively, with probability at least $1 - \delta$, we have
\begin{equation*}
    \widehat{v}(X) \leq v + \frac{4 \Delta \psi \log (1 /\delta)}{n} + 2 \sqrt{\frac{2 v \Delta \psi \log (1 /\delta)}{n}},
\end{equation*}
\begin{equation*}
    \widehat{v}(X) \leq \left( \sqrt{v} + \sqrt{\frac{2 \Delta \psi \log (1 /\delta)}{n}} \right)^2 + \frac{2 \Delta \psi \log (1 /\delta)}{n},
\end{equation*}
\begin{equation*}
    \sqrt{\widehat{v}(X)} \leq \sqrt{\left( \sqrt{v} + \sqrt{\frac{2 \Delta \psi \log (1 /\delta)}{n}} \right)^2 + \frac{2 \Delta \psi \log (1 /\delta)}{n}},
\end{equation*}
\begin{equation*}
    \sqrt{\widehat{v}(X)} \leq \sqrt{v} + 2 \sqrt{\frac{ 2\Delta \psi \log (1 /\delta)}{n}},
\end{equation*}
which finishes proving the lemma\footnote{In this proof, \citet{maurer2006concentration} is sufficient. We do not need the stronger results of \citet{boucheron2009concentration}.} for $b = 1$.

\qedsymbol

\subsection{Proof of Theorem~\ref{theorem:empirical-variance-mmd-bound}}

From Theorem~\ref{theorem:variance-aware-mmd-estimation}, 
with probability at least $1 - \delta/2$, it holds that
\begin{equation}
\label{eq:combination-first}
    \nrm{\ekme(X_1, \dots, X_n) - \kme}_{\calH_\psi}  \leq  \sqrt{2 v \frac{\log (4/ \delta)}{n}} + \frac{4 }{3}\sqrt{\overline{\psi}}\frac{\log (4/ \delta)}{n}.
\end{equation}
Invoking Lemma~\ref{lemma:concentration-square-root-v} for $b = -1$, at confidence $1 - \delta /2$ it holds that
\begin{equation}
\label{eq:combination-second}
    \sqrt{v}  \leq \sqrt{\widehat{v}(X)} +  2 \sqrt{\frac{2\Delta \psi \log (2 /\delta)}{n}} \leq  \sqrt{\widehat{v}(X)} +  2 \sqrt{\frac{2\Delta \psi \log (4 /\delta)}{n}},
\end{equation}
thus by combining \eqref{eq:combination-first} and \eqref{eq:combination-second} and a union bound yields that with probability at least $1- \delta$,
\begin{equation}
\begin{split}
    &\nrm{\ekme(X_1, \dots, X_n) - \kme}_{\calH_\psi}  \\
    \leq&  \left( \sqrt{\widehat{v}(X)} +  2 \sqrt{\frac{2\Delta \psi \log (4/ \delta)}{n}} \right)\sqrt{ \frac{2 \log (4/ \delta)}{n}} + \frac{4}{3} \sqrt{\overline{\psi}}\frac{\log (4/ \delta)}{n} \\
    =& \sqrt{2 \widehat{v}(X) \frac{\log (4/ \delta)}{n}}  + 4 \left(\sqrt{\Delta \psi} +  \frac{\sqrt{\overline{\psi}}}{3} \right) \frac{\log (4/ \delta)}{n} \\
    \leq& \sqrt{2 \widehat{v}(X) \frac{\log (4/ \delta)}{n}}  + \frac{16}{3}\sqrt{\Delta \psi} \frac{\log (4/ \delta)}{n}. \\
\end{split}
\end{equation}

\qedsymbol

\subsection{Proof of Lemma~\ref{lemma:exact-expression-second-moment-covariances}}
\begin{equation*}
\begin{split}
    \bbE\nrm{ \ekme - \kme }^2_{\calH_k}
    =&  \bbE \langle{ \frac{1}{n}\sum_{t = 1}^{n} k(X_t, \cdot) - \kme, \frac{1}{n}\sum_{s = 1}^{n} k(X_s, \cdot) - \kme \rangle}_{\calH_k}  \\
    =&  \frac{1}{n^2}\sum_{s = 1}^{n}\sum_{t = 1}^{n}\bbE  \langle{ k(X_t, \cdot) - \kme,  k(X_s, \cdot) - \kme \rangle}_{\calH_k}  \\
    =&  \frac{2}{n^2}\sum_{t = 1}^{n} \sum_{s < n}\bbE \langle{ k(X_t, \cdot) - \kme,  k(X_s, \cdot) - \kme \rangle}_{\calH_k} \\ &+  \frac{1}{n^2}\sum_{t = 1}^{n}\bbE \langle{ k(X_t, \cdot) - \kme,  k(X_t, \cdot) - \kme \rangle}_{\calH_k}  \\
    =&  \frac{2}{n^2}\sum_{t = 2}^{n} (n - t + 1) \bbE \langle{ k(X_t, \cdot) - \kme,  k(X_1, \cdot) - \kme \rangle}_{\calH_k}  +  \frac{v}{n}.
\end{split}
\end{equation*}
\qedsymbol

\subsection{Proof of Theorem~\ref{theorem:variance-aware-mmd-estimation-phi-mixing}}
An application of Jensen's inequality and Lemma~\ref{lemma:exact-expression-second-moment-covariances} immediately yield that
\begin{equation*}
    \nrm{\ekme - \kme}_{\calH_k} \leq \sqrt{ \frac{v + \Sigma_n}{n}},
\end{equation*}
where $\Sigma_n$ is a total measure of covariance in the RKHS introduced in Lemma~\ref{lemma:exact-expression-second-moment-covariances}.
Additionally, recall that we can express the norm of the deviation of $\ekme$ from its mean $\kme$ in the RKHS as a supremum of empirical processes as follows,
\begin{equation*}
    \nrm{\ekme - \kme}_{\calH_k} = \sup_{\substack{f \in \calH_k \\ \nrm{f}_{\calH_k} \leq 1}} \langle f, \ekme - \kme  \rangle_{\calH_k} = 
    \frac{1}{n}\sup_{\substack{f \in \calH_k \\ \nrm{f}_{\calH_k} \leq 1}} \sum_{t = 1}^{n} \left( f(X_t) - \bbE f \right) = \frac{1}{n} \sup_{g \in \calF} \sum_{t = 1}^{n} g(X_t),
\end{equation*}
where $\calF$ is the re-centered class,
\begin{equation*}
    \calF \eqdef \set{ g \colon \calX \to \bbR, g = f - \bbE f, f \in \calH_{k}, \nrm{f}_{\calH_k} \leq 1}. 
\end{equation*}
In particular, observe that $g \in \calF \implies \abs{g} \leq \overline{k}$. Furthermore, the variance $v$ in the RKHS can be expressed as
\begin{equation*}
\begin{split}
    v &= \bbE \nrm{k(X, \cdot) - \kme}_{\calH_k}^2 \\ &= \bbE  \left(\sup_{\substack{f \in \calH_k \\ \nrm{f}_{\calH_k} \leq 1} } \langle f, 
k(X, \cdot) - \kme \rangle_{\calH_k} \right)^2 \\
&= \bbE \sup_{\substack{f \in \calH_k \\ \nrm{f}_{\calH_k} \leq 1} } \langle f, 
k(X, \cdot) - \kme \rangle_{\calH_k}^2 \\
&= \bbE \sup_{\substack{f \in \calH_k \\ \nrm{f}_{\calH_k} \leq 1} } \left( f(X) - \bbE f \right)^2.
\end{split}
\end{equation*}
We invoke Theorem~\ref{theorem:talagrand-type-concentration-phi-mixing-processes}
for $\phi$-dependent random sequences,
with $Z = n \nrm{\ekme - \kme}_{\calH_k}$ and observing that
\begin{equation*}
    \E{V} =  \sum_{t=1}^{n} \bbE \sup_{g \in \calF} g(X_t)^2 = n v,
\end{equation*}
with probability at least $1- \delta$ it holds that
\begin{equation*}
    \nrm{\ekme - \kme}_{\calH_k} \leq \sqrt{\frac{v + \Sigma_n}{n}} + 4 \sqrt{\frac{2 v \nrm{\Gamma}_2 \log (1/\delta)}{n}} + \frac{8 \overline{k} \nrm{\Gamma}_2 \log (1/\delta)}{n}.
\end{equation*}
\qedsymbol

\subsection{Proof of Theorem~\ref{theorem:hypothesis-testing}}
First, assume that $\mathbf{H}_0$ is true, that is: there is a $\theta_0\in\Theta$ such that $\bbP=\bbP_{\theta_0}$. Then,
\begin{align*}
\bbP_{\mathbf{H}_0}\left(  \mathcal{C} \right)
& =
\bbP_{\mathbf{H}_0}\left(  T > \mathcal{B}(X_1,\dots,X_n,\alpha) \right)
\\
& = \bbP_{\mathbf{H}_0}\left(  \inf_{\theta\in\Theta}\nrm{\widehat{\mu_{\bbP}}- \mu_{\bbP_\theta}}_{\calH_k} > \mathcal{B}(X_1,\dots,X_n,\alpha) \right)
\\
& \leq \bbP_{\mathbf{H}_0}\left(  \nrm{\widehat{\mu_{\bbP}}- \mu_{\bbP_{\theta_0}}}_{\calH_k} > \mathcal{B}(X_1,\dots,X_n,\alpha) \right)
\\
& \leq \alpha
\end{align*}
thanks to $(a)$.

Now, assume that $\mathbf{H}_0$ is not true, that is $\bbP\notin \{\bbP_{\theta},\theta\in\Theta \}$. Note that under the assumption that the model is closed, we have
$$ 0< \Delta :=  \inf_{\theta\in\Theta}  \nrm{\mu_{\bbP}- \mu_{\bbP_{\theta}}}_{\calH_k}. $$
By the triangle inequality, for any $\theta$,
$$
 \nrm{\mu_{\bbP_\theta}- \widehat{\mu}_{\bbP}}_{\calH_k} \geq  \nrm{\mu_{\bbP_\theta}- \mu_{\bbP}}_{\calH_k} - \nrm{\mu_{\bbP}- \widehat{\mu}_{\bbP}}_{\calH_k}
$$
and thus, taking the infimum w.r.t $\theta$ on both sides,
$$
 T \geq  \Delta - \nrm{\mu_{\bbP}- \widehat{\mu}_{\bbP}}_{\calH_k}.
$$
We apply $(a)$ with confidence level $\alpha_n$ and obtain
$$
\bbP_{\mathbf{H}_1}\left( \nrm{\mu_{\bbP}- \widehat{\mu}_{\bbP}}_{\calH_k} \geq \mathcal{B}(X_1,\dots,X_n,\alpha_n) \right) \leq \alpha_n
$$
and thus
$$
\bbP_{\mathbf{H}_1}\left( T \leq \Delta - \mathcal{B}(X_1,\dots,X_n,\alpha_n) \right) \leq \alpha_n.
$$
As $ \mathcal{B}(X_1,\dots,X_n,\alpha_n)\rightarrow 0$ when $n\rightarrow\infty$ thanks to $(c)$, there is a $N$ large enough such that, for any $n\geq N$, $\mathcal{B}(X_1,\dots,X_n,\alpha_n)\leq \Delta/2$, and thus
$$
\bbP_{\mathbf{H}_1}\left( T \leq \frac{\Delta}{2} \right) \leq \alpha_n.
$$
Moreover, as we also have $\mathcal{B}(X_1,\dots,X_n,\alpha_n)\rightarrow 0$ when $n\rightarrow\infty$ thanks to $(b)$, there is a $N'$ large enough such that, for any $n\geq N'$, $\mathcal{B}(n,\alpha)\leq \Delta/2$ and thus
$$
\bbP_{\mathbf{H}_1}\left( T \leq \mathcal{B}(X_1,\dots,X_n,\alpha_n)\right) \leq \alpha_n.
$$
Finally,
$$
\bbP_{\mathbf{H}_1}\left( \mathcal{C} \right) = 1 - \bbP_{\mathbf{H}_1}\left( T \leq \mathcal{B}(X_1,\dots,X_n,\alpha_n) \right) \geq 1-\alpha_n
$$
as soon as $n\geq \max(N,N')$.

\qedsymbol

\subsection{Proof of Theorem~\ref{theorem:hypothesis-testing_2}}
First, under $\mathbf{H}_0$, we have $\bbP_X=\bbP_Y$ and thus
\begin{align*}
\bbP_{\mathbf{H}_0}\left(  \mathcal{C}_2 \right)
& =
\bbP_{\mathbf{H}_0}\left(  T_2 > \mathcal{B}(X_1,\dots,X_n,\alpha/2) + \mathcal{B}(Y_1,\dots,Y_m,\alpha/2)\right)
\\
& = \bbP_{\mathbf{H}_0}\Biggl( \nrm{\widehat{\mu_{\bbP}}(X_1,\dots,X_n)- \widehat{\mu_{\bbP}}(Y_1,\dots,Y_m)}_{\calH_k} 
\\
&  \hspace{3cm} > \mathcal{B}(X_1,\dots,X_n,\alpha/2) + \mathcal{B}(Y_1,\dots,Y_m,\alpha/2)\Biggr)
\\
& = \bbP_{\mathbf{H}_0}\Biggl( \nrm{\widehat{\mu_{\bbP}}(X_1,\dots,X_n) - \mu_{\bbP_X} + \mu_{\bbP_X} - \widehat{\mu_{\bbP}}(Y_1,\dots,Y_m)}_{\calH_k} 
  \\
  & \hspace{3cm} > \mathcal{B}(X_1,\dots,X_n,\alpha/2) + \mathcal{B}(Y_1,\dots,Y_m,\alpha/2)  \Biggr)
\\
& \leq \bbP_{\mathbf{H}_0}\left( \nrm{\widehat{\mu_{\bbP}}(X_1,\dots,X_n) - \mu_{\bbP_X} }_{\calH_k} > \mathcal{B}(X_1,\dots,X_n,\alpha/2)  \right) \\
& \quad \quad \quad + \bbP_{\mathbf{H}_0}\left( \nrm{\widehat{\mu_{\bbP}}(Y_1,\dots,Y_m) - \mu_{\bbP_X} }_{\calH_k} > \mathcal{B}(Y_1,\dots,Y_m,\alpha/2)  \right)
\\
& \leq \alpha/2 + \alpha/2 = \alpha.
\end{align*}
thanks to $(a)$.

Now, assume that $\mathbf{H}_0$ is not true, that is $\bbP_X\neq \bbP_Y$, and put
$$ \Delta :=    \nrm{\mu_{\bbP_X}- \mu_{\bbP_Y}}_{\calH_k}>0. $$
By the triangle inequality,
\begin{multline*}
 T_2 = \nrm{ \widehat{\mu_{\bbP}}(X_1,\dots,X_n)- \widehat{\mu_{\bbP}}(Y_1,\dots,Y_m) }_{\calH_k} \\
 \geq  \Delta - \nrm{\widehat{\mu_{\bbP}}(X_1,\dots,X_n) - \mu_{\bbP_X} }_{\calH_k} - \nrm{\widehat{\mu_{\bbP}}(Y_1,\dots,Y_n) - \mu_{\bbP_Y} }_{\calH_k}.
\end{multline*}
We apply $(a)$ with respectice confidence levels $\alpha_n$ and $\alpha_m$ to get
$$
\bbP_{\mathbf{H}_1}\left( \nrm{\mu_{\bbP_X}- \widehat{\mu}_{\bbP}(X_1,\dots,X_n)}_{\calH_k} \geq \mathcal{B}(X_1,\dots,X_n,\alpha_n) \right) \leq \alpha_n
$$
and
$$\bbP_{\mathbf{H}_1}\left( \nrm{\mu_{\bbP_Y}- \widehat{\mu}_{\bbP}(Y_1,\dots,Y_m)}_{\calH_k} \geq \mathcal{B}(Y_1,\dots,Y_m,\alpha_m) \right) \leq \alpha_m.
$$
Thus,
$$
\bbP_{\mathbf{H}_1}\left( T_2 \leq \Delta - \mathcal{B}(X_1,\dots,X_n,\alpha_n) - \mathcal{B}(Y_1,\dots,Y_m,\alpha_m) \right) \leq \alpha_n + \alpha_m.
$$
Thanks to $(c)$, there is a $N$ large enough such that, as soon as both $n,m\geq N$, both  $\mathcal{B}(X_1,\dots,X_n,\alpha_n)\leq \Delta/4 $ and $\mathcal{B}(Y_1,\dots,Y_m,\alpha_m)\leq \Delta/4$ hold, and thus
$$
\bbP_{\mathbf{H}_1}\left( T \leq \frac{\Delta}{2} \right) \leq \alpha_n + \alpha_m.
$$
Moreover, thanks to $(b)$, there is $N'$ large enough such that, when $n,m\geq N'$, both $\mathcal{B}(X_1,\dots,X_n,\alpha/2)\leq \Delta/4$ and $\mathcal{B}(Y_1,\dots,Y_m,\alpha/2)\leq \Delta/4$ hold, and thus
$$
\bbP_{\mathbf{H}_1}\left( T_2 \leq  \mathcal{B}(X_1,\dots,X_n,\alpha/2) + \mathcal{B}(Y_1,\dots,Y_m,\alpha/2) \right) \leq \alpha_n + \alpha_m.
$$
Finally,
$$
\bbP_{\mathbf{H}_1}\left( \mathcal{C}_2 \right) = 1 - \bbP_{\mathbf{H}_1}\left( T_2 \leq \mathcal{B}(X_1,\dots,X_n,\alpha/2) + \mathcal{B}(Y_1,\dots,Y_m,\alpha/2) \right) \geq 1-\alpha_n-\alpha_m
$$
as soon as $n\geq \max(N,N')$, and thus
$$
\bbP_{\mathbf{H}_1}\left( \mathcal{C}_2 \right) \geq 1-\alpha_n-\alpha_m \xrightarrow[n,m\rightarrow\infty]{} 1.
$$
\qedsymbol

\subsection{Proof of Lemma~\ref{lemma:variance-huber-contamination}}
\begin{equation*}
\begin{split}
    v(\bbP) &= \bbE_{X \sim \bbP} \nrm{k(X, \cdot) - \kme}_{\calH_k}^2 \\
    &= \overline{\psi} - \int \int \psi(x' - x) d(\bbP \times \bbP)(x,x') \\
    &= \overline{\psi} - \int \int \psi(x' - x) d\bbP(x)d\bbP(x') \\
    &= \overline{\psi} - \int \int \psi(x' - x) d((1 - \xi)\bbP_{\theta_0}(x) + \xi \bbH(x))d((1 - \xi)\bbP_{\theta_0}(x') + \xi \bbH(x')) \\
    &= \overline{\psi} - (1 - \xi)^2 \int \int \psi(x' - x) d\bbP_{\theta_0}(x)d\bbP_{\theta_0}(x') \\
    &\qquad- 2 \xi(1 - \xi) \int \int \psi(x' - x) d\bbP_{\theta_0}(x)d\bbH(x') \\
    &\qquad- \xi^2 \int \int \psi(x' - x) d\bbH(x)d\bbH(x') \\
    &\leq \overline{\psi} - (1 - \xi)^2 (\overline{\psi} - v) - 2 \xi(1 - \xi) \underline{\psi} - \xi^2 \underline{\psi}\\
    &= \overline{\psi} - (1 - 2\xi + \xi^2) (\overline{\psi} - v) - \xi (2 - \xi) \underline{\psi}\\
    &= \overline{\psi} - (\overline{\psi} - v) + 2 \xi (\overline{\psi} - v) - \xi^2 (\overline{\psi} - v)  - \xi (2 - \xi) \underline{\psi}\\
    &\leq (1 - 2 \xi) v + 2 \xi \overline{\psi} - \xi (2 - \xi) \underline{\psi},\\
    &= (1 - 2 \xi) v + 2 \xi \Delta \psi + \xi^2 \underline{\psi},\\
\end{split}
\end{equation*}
and the lemma follows by disjunction of cases.
\qedsymbol

\subsection*{Acknowledgments}
The authors declare that there is no conflict of interest. Part of this research was conducted when GW was supported by the Special Postdoctoral Researcher Program (SPDR) of RIKEN. Part of the research in this paper was done when PA was working at RIKEN AIP. We thank Bastien Dussap for pointing out a missing factor two in an earlier version of this manuscript. We thank the anonymous reviewers for their useful comments. In particular, the comparison to the approach based on Bernstein's inequality for U-statistics in Section~\ref{section:two-means} was suggested by Reviewer 1 and the analysis for time-dependent processes in Section 5 was suggested by Reviewer 2. All the remaining mistakes are ours.

\bibliography{bibliography}
\bibliographystyle{abbrvnat}

\appendix

\section{Extension to Non Translation Invariant Kernels}
\label{section:non-translation-invariant-kernels}

Let us define $\diag k \colon \calX \to \R$ by
$\diag k(x) = k(x,x)$, and use the shorthands $\underline{\diag k}, \overline{\diag k}$ and $\Delta \diag k$ introduced in \eqref{definition:shorthands} accordingly. We show that the properties obtained in Section~\ref{section:variance-empirical} can be extended to non TI kernel when $\Delta \diag k$ is controlled from above.
\begin{lemma}
\label{lemma:general-weakly-self-bounding}
Let $k$ be a characteristic reproducing kernel.
The function
\begin{equation*}
    \calX^n \to \R, x \mapsto \frac{n}{\Delta \diag k + 2 \Delta k} \widehat{v}(x),
\end{equation*}
is weakly $$\left(2, 2 n \frac{ \Delta \diag k}{\Delta \diag k + 2 \Delta k} \right)\text{-self-bounding},$$ and has bounded differences in the sense of Definition~\ref{definition:weakly-self-bounding-and-bounded-differences-function}.
\end{lemma}

\begin{lemma}
\label{lemma:general-concentration-square-root-v}
For $b \in \set{-1,1}$, with probability at least $1- \delta$,
\begin{equation}
\label{eq:general-square-root-confidence-interval-for-v}
    b \left[ \sqrt{\widehat{v}(X) + \Delta \diag k} - \sqrt{v + \Delta \diag k } \right] \leq 2 \sqrt{\frac{(\Delta \diag k + 2 \Delta k) \log (1 /\delta)}{n}}.
\end{equation}
\end{lemma}

Attempting to recover concentration of $\sqrt{\widehat{v}}$ around $\sqrt{v}$ with self-boundedness leads to an additional term in $\bigO\left( n^{-1/4} \right)$ in \eqref{eq:general-square-root-confidence-interval-for-v} when $\Delta \diag k \neq 0$. 
We thus settle for concentration of $\sqrt{\widehat{v} + \Delta \diag k}$ around $\sqrt{v + \Delta \diag k}$ instead.

\begin{theorem}[Confidence interval with empirical variance for general kernel]

\label{theorem:general-empirical-variance-mmd-bound}
Let $n \in \N$, $X_1, \dots, X_n \sim \bbP$, let $k$ be a characteristic reproducing kernel defined from a positive definite function $\psi$ [see \eqref{definition:translation-invariant}].
Then with probability at least $1- \delta$, it holds that
\begin{equation*}
\begin{split}
    \nrm{\ekme - \kme}_{\calH_{k}} &\leq \widehat{\calB}_{k, \delta}(X_1, \dots, X_n), \\
    \end{split}
    \end{equation*}
    with 
    \begin{equation*}
\begin{split}
    \widehat{\calB}_{k, \delta}(X_1, \dots, X_n) = \widehat{\calB}_{\delta} \eqdef& \sqrt{2 (\widehat{v}(X_1, \dots, X_n) + \Delta \diag k) \frac{\log (4/ \delta)}{n}} \\  &+ \left(\frac{16}{3} \sqrt{\Delta k} + 2\sqrt{2} \sqrt{\Delta \diag k}\right)  \frac{ \log (4 /\delta)}{n},
\end{split}
\end{equation*}
and where $\widehat{v}$ is the empirical variance proxy defined in \eqref{eq:general-empirical-variance-proxy}.
\end{theorem}

\begin{remark}
In particular, for a TI kernel $\Delta \diag k = 0$, recovering the results of the previous section. It is currently unclear whether the term in $\Delta \diag k$ is necessary or if it is an artifact of our proof.
\end{remark}

\section{Additional Proofs}

In this section, we present supplementary technical proofs that were postponed in previous sections of the paper

\subsection{Proof of Lemma~\ref{lemma:epanechnikov-self-bounding}}

For some $t_0$, we let $x^{(t_0)} = (x_1, \dots, x_{t_0 - 1}, x'_{t_0}, x_{t_0 +1}, \dots, x_n)$ where $x_{t_0}$ was replaced with $x'_{t_0}$. 
\begin{equation}
\label{eq:covariance-matrix-trace-bounded-differences}
\begin{split}
    &\frac{1}{d} \Trace  \widehat{\Sigma}(x) - \frac{1}{d} \Trace  \widehat{\Sigma}\left(x^{(t_0)}\right) \\
    =& \frac{1}{d}\sum_{i = 1}^{d}
    \frac{1}{2 n (n-1)} \left(\sum_{s = 1}^{n}\sum_{t = 1}^{n} \left(x_t^i - x_s^i\right)^2 -  \sum_{s = 1}^{n}\sum_{t = 1}^{n} \left(x_t^{(t_0)i} - x_s^{(t_0)i}\right)^2 \right) \\
    \leq& \frac{1}{d} \sum_{i = 1}^d \left( \frac{1}{n-1} \sum_{s = 1}^n \left( x_{t_0}^i - x_s^i\right)^2 \right) \leq 1,
\end{split}
\end{equation}
thus 
\begin{equation*}
\begin{split}
    \frac{1}{d} \Trace  \widehat{\Sigma}(x) - \inf_{x'_{t_0} \in \calX} \frac{1}{d} \Trace  \widehat{\Sigma}\left(x^{(t_0)}\right) \leq 1.
\end{split}
\end{equation*}
Furthermore,
\begin{equation*}
\begin{split}
    \sum_{t_0 = 1}^{n} \left( \frac{1}{d} \Trace  \widehat{\Sigma}(x) - \frac{1}{d} \Trace  \widehat{\Sigma}\left(x^{(t_0)}\right) \right)^2 &\stackrel{(i)}{\leq}\sum_{t_0 = 1}^{n} \left(\frac{1}{d} \sum_{i = 1}^d \left( \frac{1}{n-1} \sum_{s = 1}^n \left( x_{t_0}^i - x_s^i\right)^2 \right) \right)^2 \\
    &\stackrel{(ii)}{\leq} \frac{1}{d} \sum_{i = 1}^d \sum_{t_0 = 1}^{n} \left( \frac{1}{n-1} \sum_{s = 1}^n \left( x_{t_0}^i - x_s^i\right)^2  \right)^2 \\
    &\stackrel{(iii)}{\leq} \frac{1}{d} \sum_{i = 1}^d \frac{n}{n-1} \widehat{v}^i(X) \\
    &= \frac{n}{n-1} \frac{1}{d} \Trace  \widehat{\Sigma}(X),
\end{split}
\end{equation*}
where $(i)$ follows from 
\eqref{eq:covariance-matrix-trace-bounded-differences}, $(ii)$ is Jensen's inequality, and $(iii)$ stems from \citet[Corollary~9]{maurer2009empirical}.

\qedsymbol

\subsection{Proof of Lemma~\ref{lemma:general-empirical-proxy-unbiased}}
We verify that $\widehat{v}$ is an unbiased estimator for $v$.
\begin{equation*}
\begin{split}
    &\E[X_1, \dots, X_n \sim \bbP]{\widehat{v}(X_1, \dots, X_n)} \\  =& \E[X_1, \dots, X_n \sim \bbP]{\frac{1}{n-1} \sum_{t = 1}^{n} \left( k(X_t, X_t) - \frac{1}{n}\sum_{s=1}^{n} k(X_t, X_s) \right)} \\
    =& \frac{1}{n-1} \sum_{t = 1}^{n} \left( \E[X \sim \bbP]{k(X, X)} - \frac{1}{n}\sum_{s=1}^{n} \E[X_t, X_s \sim \bbP]{k(X_t, X_s)} \right) \\
    =& \frac{n}{n-1}  \E[X \sim \bbP]{k(X, X)} - \frac{1}{n(n-1)}\sum_{t=1}^{n} \sum_{s=1}^{n} \E[X_t, X_s \sim \bbP]{k(X_t, X_s)}  \\
    =& \frac{n}{n-1}  \E[X \sim \bbP]{k(X, X)} - \frac{1}{(n-1)} \E[X \sim \bbP]{k(X, X)} \\ 
    & \qquad - \frac{n^2 - n}{n(n-1)} \E[X,X' \sim \bbP]{k(X, X')} \\
    =& \E[X \sim \bbP]{k(X, X)} - \E[X,X' \sim \bbP]{k(X, X')} \\
    =& v.
\end{split}
\end{equation*}

\qedsymbol

\subsection{Proof of Theorem~\ref{theorem:variance-aware-mmd-estimation-beta-mixing}}
For $t \in \bbN$, we write $Z_t \eqdef k(X_t, \cdot) - \kme$, thus
\begin{equation*}
    \nrm{\ekme(X_1, \dots, X_n) - \kme}_{\calH_k} = \frac{1}{n} \nrm{\sum_{t=1}^{n}Z_t}_{\calH_k}.
\end{equation*}
We assume that $n = 2Bs$, where $B, s \in \bbN$ will be determined later. For $b \in [B]$, we denote
\begin{equation*}
\begin{split}
    Z^{[2b]} &\eqdef  Z^{[2b]}_1, \dots, Z^{[2b]}_s  \eqdef Z_{(2b - 1)s + 1}, \dots, Z_{2bs}, \\
    Z^{[2b - 1]} &\eqdef  Z^{[2b - 1]}_1, \dots, Z^{[2b - 1]}_s  \eqdef  Z_{(2b - 2)s + 1}, \dots, Z_{(2b - 1)s} , \\
\end{split}
\end{equation*}
and we decompose the random process into blocks
\begin{equation*}
   Z^{[1]}, Z^{[2]}, \dots, Z^{[2B]}. 
\end{equation*}
From sub-additivity of the norm and the blocking method described in \citet[Corollary~2.7]{yu1994rates}, it holds that
\begin{equation*}
\begin{split}
\PR{\nrm{\sum_{t=1}^{n}Z_t}_{\calH_k} > n \eps} &\leq \sum_{\sigma \in \set{0,1}} \PR{ \nrm{\sum_{b = 1}^{B} \left(\sum_{t = 1}^s Z_t^{[2b - \sigma]}\right)}_{\calH_k} > n \eps/2 } \\ 
&\leq \sum_{\sigma \in \set{0,1}} \PR{ \nrm{  \sum_{b = 1}^{B} \left( \sum_{t = 1}^s \widetilde{Z}_t^{[2b - \sigma]}\right)}_{\calH_k} > n \eps/2 } + 2 \left(B - 1\right) \beta(s) \\
\end{split}
\end{equation*}
where for $b \in [2]$, the blocks $\widetilde{Z}^{[2b]} = \widetilde{Z}^{[2b]}_1, \dots, \widetilde{Z}^{[2b]}_s$ are mutually independent and a similar fact holds for the blocks $\widetilde{Z}^{[2b - 1]} = \widetilde{Z}^{[2b - 1]}_1, \dots, \widetilde{Z}^{[2b - 1]}_s$.
As a result, for $b \in [B]$, the $\sum_{t = 1}^s \widetilde{Z}_t^{[2b - \sigma]}$ are independent, centered and valued in $\calH_k$, thus the problem above has been reduced to controlling the norm of a sum of $B$ iid random vectors in the Hilbert space $\calH_k$.
Let us fix $\sigma \in \set{0, 1}$.
Similar to our proof in the iid setting, our goal is to invoke Theorem~\ref{theorem:bernstein-inequality-hilbert-space}
\citep[Theorem~3.3.4]{yurinsky1995sums}.
In order to do that, we need to control from above the $\sum_{b = 1}^{B} \bbE \nrm{ \sum_{t = 1}^s \widetilde{Z}_t^{[2b - \sigma]}}_{\calH_k}^p$,
for any $p \geq 3$.
For $p=2$, by linearity of the expectation, we have
\begin{equation*}
\begin{split}
    \sum_{b = 1}^{B} \bbE \nrm{ \sum_{t = 1}^s \widetilde{Z}_t^{[2b - \sigma]}}_{\calH_k}^2 
    &= \sum_{b = 1}^{B} \sum_{t = 1}^{s} \sum_{r = 1}^{s} \bbE \langle  \widetilde{Z}_t^{[2b - \sigma]},  \widetilde{Z}_r^{[2b - \sigma]} \rangle_{\calH_k}. \\
\end{split}
\end{equation*}
By stationarity, we can write more simply
\begin{equation*}
\begin{split}
    \sum_{b = 1}^{B} \bbE \nrm{ \sum_{t = 1}^s \widetilde{Z}_t^{[2b - \sigma]}}_{\calH_k}^2
    =& \sum_{b = 1}^{B} \sum_{t = 1}^{s} \sum_{r = 1}^{s} \bbE \langle  k(X_{(2b - \sigma - 1)s + t}, \cdot) - \kme, k(X_{(2b - \sigma - 1)s + r},\cdot) - \kme \rangle_{\calH_k} \\
    =& B \sum_{t = 1}^{s} \sum_{r = 1}^{s} \bbE \langle  k(X_{t - r + 1}, \cdot) - \kme, k(X_{1},\cdot) - \kme \rangle_{\calH_k} \\
    =& B \Bigg( \sum_{t = 1}^{s} \langle  k(X_{1}, \cdot) - \kme, k(X_{1},\cdot) - \kme \rangle_{\calH_k} \\ &+  2\sum_{t = 2}^{s} (s - t + 1) \bbE \langle  k(X_{t}, \cdot) - \kme, k(X_{1},\cdot) - \kme \rangle_{\calH_k} \Bigg). \\
\end{split}
\end{equation*}
It is then natural to consider the following covariance coefficients in the RKHS introduced by \cite{cherief2022finite},
\begin{equation*}
    \rho_{t} \eqdef \bbE \langle  k(X_{t}, \cdot) - \kme, k(X_{1},\cdot) - \kme \rangle_{\calH_k},
\end{equation*}
and in particular, observe that we recover $\rho_1 = v$. It follows that
\begin{equation*}
\begin{split}
    \sum_{b = 1}^{B} \bbE \nrm{ \sum_{t = 1}^s \widetilde{Z}_t^{[2b - \sigma]}}_{\calH_k}^2
    &\leq \frac{n}{2} \left( v + \frac{2}{s} \sum_{t=2}^{s} (s - t + 1) \rho_{t} \right) = \frac{n}{2} \left(v + \Sigma_s \right),
\end{split}
\end{equation*}
where the mixing coefficient $\Sigma_s$ was introduced in Lemma~\ref{lemma:exact-expression-second-moment-covariances}.
In particular, $\Sigma_s$ vanishes for iid processes.
For $p \geq 3$, by sub-additivity of the norm, 
\begin{equation*}
\begin{split}
    \nrm{ \sum_{t = 1}^s \widetilde{Z}_t^{[2b - \sigma]}}_{\calH_k} &\leq \sum_{t = 1}^s \nrm{  \widetilde{Z}_t^{[2b - \sigma]}}_{\calH_k} \leq 2 s \sqrt{\overline{k}}.
\end{split}
\end{equation*}
Therefore, invoking Theorem~\ref{theorem:bernstein-inequality-hilbert-space} with $G^2 = \frac{n}{2} (v + \Sigma_s)$ and $H = 2 s \sqrt{\overline{k}}/3$,
we obtain for $\sigma \in \set{0,1}$ and for any $\eps > 0$ that,
\begin{equation*}
\begin{split}
    \PR{ \nrm{\sum_{b = 1}^{B} \left( \sum_{t = 1}^s \widetilde{Z}_t^{[2b - \sigma]}\right)}_{\calH_k} > \frac{n \eps}{2} } &\leq 2\exp \left(- \frac{((\eps/2) n)^2}{2(n (v + \Sigma_s) /2 + ((\eps/2) n) 2 s \sqrt{\overline{k}}/3)} \right)\\
    &= 2\exp \left(- \frac{n\eps^2}{4\left(v + \Sigma_s + 2 \eps s \sqrt{\overline{k}}/3\right)} \right).
\end{split}
\end{equation*}
It remains to control the term involving the $\beta$-mixing coefficient by suitable choosing the blocking size. For
\begin{equation*}
    \tau^{\beta}_{n, \delta} \eqdef \argmin_{s \in \bbN} \set{ s \geq \betamix \left( \frac{\delta}{6(n/(2s) - 1)} \right)},
\end{equation*}
it holds that $2 \left(B - 1\right) \beta(s) \leq \delta/3$, which concludes the proof.

\qedsymbol

\subsection{Proof of Lemma~\ref{lemma:general-weakly-self-bounding}}
We will rely on the following property for general kernels.

\begin{lemma}
\label{lemma:general-empirical-proxy-bounded-differences}
Let $x = (x_1, \dots, x_{t_0}, \dots, x_n) \in \calX^n$ and $x_{t_0}' \in \calX$. Writing 
$$x^{(t)} = (x_1, \dots, x_{{t_0}-1}, x_{t_0}', x_{{t_0}+1}, \dots, x_n),$$ it holds that
\begin{equation*}
\begin{split}
    \widehat{v}(x) - \widehat{v}\left(x^{({t_0})}\right) =& \frac{1}{n} \left( k\left(x_{t_0}, x_{t_0}\right) - k\left(x'_{t_0}, x'_{t_0}\right) \right) \\ 
    &+ \frac{2}{(n-1)n} \sum_{t \in [n], t \neq {t_0}} \left( k\left(x'_{t_0}, x_t\right) - k\left(x_{t_0}, x_t\right) \right).
\end{split}
\end{equation*}
\end{lemma}

From Lemma~\ref{lemma:general-empirical-proxy-bounded-differences}, if follows that
\begin{equation*}
\begin{split}
    \widehat{v}(x) - \inf_{x_{t_0}' \in \calX} \widehat{v}\left( x^{(t_0)} \right) \stackrel{(\dagger)}{\leq}& \frac{1}{n} \left( k\left(x_{t_0}, x_{t_0}\right) - \underline{\diag k} \right) \\ 
    &+ \frac{2}{(n-1)n} \sum_{t \in [n], t \neq t_0} \left( \overline{k} - k\left(x_{t_0}, x_t\right) \right) \\
    \stackrel{(\ddagger)}{\leq}& \frac{\Delta \diag k + 2 \Delta k}{n},
\end{split}
\end{equation*}
and the bounded differences property follows directly from $(\ddagger)$.
Furthermore, as a result of $(\dagger)$,
\begin{equation*}
\begin{split}
    &\sum_{t_0 = 1}^{n} \left(\widehat{v}(x) - \inf_{x_{t_0}' \in \calX} \widehat{v}\left( x^{(t_0)} \right) \right)^2 \\
    \leq& \sum_{t_0 = 1}^{n} \left(\frac{1}{n} \left( k\left(x_{t_0}, x_{t_0}\right) - \underline{\diag k} \right) + \frac{2}{(n-1)n} \sum_{t \in [n], t \neq t_0} \left( \overline{k} - k\left(x_{t_0}, x_t\right) \right) \right)^2 \\
    =& \frac{1}{n^2}\sum_{t_0 = 1}^{n}  \left( \frac{1}{(n-1)} \sum_{t \in [n], t \neq t_0} k\left(x_{t_0}, x_{t_0}\right) - \underline{\diag k} +  2\overline{k} - 2k\left(x_{t_0}, x_t\right) \right)^2 \\
    \leq& \frac{1}{n^2}\sum_{t_0 = 1}^{n} \frac{1}{(n-1)} \sum_{t \in [n], t \neq t_0} \left( k\left(x_{t_0}, x_{t_0}\right) - \underline{\diag k} +  2\overline{k} - 2k\left(x_{t_0}, x_t\right) \right)^2 \\
    \leq& \frac{\Delta \diag k + 2 \Delta k}{n^2}\sum_{t_0 = 1}^{n} \frac{1}{(n-1)} \sum_{t \in [n], t \neq t_0} \left( k\left(x_{t_0}, x_{t_0}\right) - \underline{\diag k} +  2\overline{k} - 2k\left(x_{t_0}, x_t\right) \right) \\
    \leq& \frac{\Delta \diag k + 2 \Delta k}{n^2}\sum_{t_0 = 1}^{n} \frac{1}{(n-1)} \sum_{t \in [n], t \neq t_0} \left( 2k\left(x_{t_0}, x_{t_0}\right) - 2\underline{\diag k} +  2\overline{k} - 2k\left(x_{t_0}, x_t\right) \right) \\
    =& 2 \frac{\Delta \diag k + 2 \Delta k}{n^2}\sum_{t_0 = 1}^{n} \frac{1}{(n-1)} \sum_{t \in [n], t \neq t_0} \left( k\left(x_{t_0}, x_{t_0}\right) - k\left(x_{t_0}, x_t\right) \right) \\
    & + 2 \frac{\Delta \diag k + 2 \Delta k}{n}\left( 
    \overline{k} - \underline{\diag k}  \right) \\
    \stackrel{(\star)}{=}& 2 \frac{\Delta \diag k + 2 \Delta k}{n^2}\sum_{t_0 = 1}^{n} \frac{1}{(n-1)} \sum_{t \in [n]} \left( k\left(x_{t_0}, x_{t_0}\right) - k\left(x_{t_0}, x_t\right) \right) \\ &+ 2 \frac{\Delta \diag k + 2 \Delta k}{n} \Delta \diag k \\
    =& 2 \frac{\Delta \diag k + 2 \Delta k}{n}\widehat{v}(x) + 2 \frac{\Delta \diag k + 2 \Delta k}{n} \Delta \diag k, \\
\end{split}
\end{equation*}
where for $(\star)$ we relied on $\overline{k} = \overline{\diag k}$ for a positive definite kernel.

\qedsymbol

\subsection{Proof of Lemma~\ref{lemma:general-concentration-square-root-v}}
Lemma~\ref{lemma:general-weakly-self-bounding} established that the function
\begin{equation*}
    \calX^n \to \R, x \mapsto \frac{n}{\Delta \diag k + 2 \Delta k} \widehat{v}(x),
\end{equation*}
is weakly self-bounding and has the bounded differences property.
From an application of \citet[Theorem~1]{boucheron2009concentration},
\begin{equation*}
\begin{split}
    \PR{\widehat{v} - v > t} &\leq \exp \left( - \frac{n t^2}{4(\Delta \diag k + 2 \Delta k)(v + \Delta \diag k + t/2)} \right), \\
    \PR{v - \widehat{v} > t} &\leq \exp \left( - \frac{n t^2}{4(\Delta \diag k + 2 \Delta k)(v + \Delta \diag k)} \right).
\end{split}
\end{equation*}
The claim follows from
rewriting
\begin{equation*}
    \widehat{v} - v = [\widehat{v} + \Delta \diag k] - [v + \Delta \diag k],
\end{equation*}
and solving quadratic inequalities as in Lemma~\ref{lemma:concentration-square-root-v}.
\qedsymbol

\subsection{Proof of Theorem~\ref{theorem:general-empirical-variance-mmd-bound}}
From Theorem~\ref{theorem:variance-aware-mmd-estimation}, and since $\Delta \diag k > 0$,
with probability at least $1 - \delta/2$, it holds that
\begin{equation}
\label{eq:general-combination-first}
    \nrm{\ekme(X_1, \dots, X_n) - \kme}_{\calH_k}  \leq  \sqrt{2 (v + \Delta \diag k) \frac{\log (4/ \delta)}{n}} + \frac{4 }{3}\sqrt{\overline{k}}\frac{\log (4/ \delta)}{n}.
\end{equation}
By Lemma~\ref{lemma:general-concentration-square-root-v} for $b = -1$, at confidence $1 - \delta /2$ it holds that
\begin{equation}
\label{eq:general-combination-second}
    \sqrt{v + \Delta \diag k}  \leq \sqrt{\widehat{v} + \Delta \diag k} +  2 \sqrt{\frac{(\Delta \diag k + 2\Delta k) \log (4/ \delta)}{n}},
\end{equation}
thus by combining \eqref{eq:general-combination-first} and \eqref{eq:general-combination-second} and a union bound yields that with probability at least $1- \delta$,
\begin{equation}
\begin{split}
    &\nrm{\ekme - \kme}_{\calH_k}  \\
    \leq&  \left( \sqrt{\widehat{v} + \Delta \diag k} +  2 \sqrt{\frac{(\Delta \diag k + 2\Delta k) \log (4/ \delta)}{n}} \right)\sqrt{ \frac{2 \log (4/ \delta)}{n}} + \frac{4}{3} \sqrt{\overline{k}}\frac{\log (4/ \delta)}{n} \\
    =& \sqrt{2 (\widehat{v} + \Delta \diag k) \frac{\log (4/ \delta)}{n}}  +  \left(2 \sqrt{2} \sqrt{\Delta \diag k + 2\Delta k} +  \frac{4}{3}\sqrt{\overline{k}} \right) \frac{ \log (4/ \delta)}{n} \\
    \leq& \sqrt{2 (\widehat{v} + \Delta \diag k) \frac{\log (4/ \delta)}{n}}  + \left(\frac{16}{3} \sqrt{\Delta k} + 2\sqrt{2} \sqrt{\Delta \diag k}\right)  \frac{ \log (4/ \delta)}{n},
\end{split}
\end{equation}
where the last inequality is by sub-additivity of the square root.
\qedsymbol

\subsection{Proof of Lemma~\ref{lemma:general-empirical-proxy-bounded-differences}}

\begin{equation*}
\begin{split}
    &\widehat{v}(x) - \widehat{v}\left(x^{(t)}\right) \\
    =& \frac{1}{n-1} \sum_{r = 1}^{n} \left( k\left(x_r, x_r\right) - \frac{1}{n}\sum_{s=1}^{n} k\left(x_r, x_s\right) \right) \\ &-\frac{1}{n-1} \sum_{r = 1}^{n} \left( k\left(x^{(t)}_r, x^{(t)}_r\right) - \frac{1}{n}\sum_{s=1}^{n} k\left(x^{(t)}_r, x^{(t)}_s\right) \right) \\
    =& \frac{1}{n-1} \sum_{r = 1}^{n} \left( k\left(x_r, x_r\right) -   k\left(x^{(t)}_r, x^{(t)}_r\right) \right) \\ &-\frac{1}{(n-1)n} \sum_{(r, s) \in [n]^2} \left( k\left(x_r, x_s\right) - k\left(x^{(t)}_r, x^{(t)}_s\right) \right)
\end{split}
\end{equation*}
thus,
\begin{equation*}
\begin{split}
&\widehat{v}(x) - \widehat{v}\left(x^{(t)}\right) \\
        =& \frac{1}{n-1} \left( k\left(x_t, x_t\right) -   k\left(x'_t, x'_t\right) \right) -\frac{1}{(n-1)n} \sum_{(r, s) \in [n]^2, r = s} \left( k\left(x_r, x_s\right) - k\left(x^{(t)}_r, x^{(t)}_s\right) \right) \\
    &-\frac{1}{(n-1)n} \sum_{(r, s) \in [n]^2, r \neq s} \left( k\left(x_r, x_s\right) - k\left(x^{(t)}_r, x^{(t)}_s\right) \right), \\
    =& \frac{1}{n-1} \left( k\left(x_t, x_t\right) -   k\left(x'_t, x'_t\right) \right) -\frac{1}{(n-1)n} \sum_{r=1}^{n} \left( k\left(x_r, x_r\right) - k\left(x^{(t)}_r, x^{(t)}_r\right) \right) \\
    &-\frac{1}{(n-1)n} \sum_{(r, s) \in [n]^2, r \neq s, r = t} \left( k\left(x_r, x_s\right) - k\left(x^{(t)}_r, x^{(t)}_s\right) \right) \\
    &- \frac{1}{(n-1)n} \sum_{(r, s) \in [n]^2, r \neq s, s = t} \left( k\left(x_r, x_s\right) - k\left(x^{(t)}_r, x^{(t)}_s\right) \right) \\
    =& \frac{1}{n} \left( k\left(x_t, x_t\right) - k\left(x'_t, x'_t\right) \right) \\
    &-\frac{1}{(n-1)n} \sum_{s \in [n], s \neq t} \left( k\left(x_t, x_s\right) - k\left(x'_t, x_s\right) \right) \\ &- \frac{1}{(n-1)n} \sum_{s \in [n], s \neq t} \left( k\left(x_t, x_s\right) - k\left(x'_t, x_s\right) \right) \\
    =& \frac{1}{n} \left( k\left(x_t, x_t\right) - k\left(x'_t, x'_t\right) \right) + \frac{2}{(n-1)n} \sum_{s \in [n], s \neq t} \left( k\left(x'_t, x_s\right) - k\left(x_t, x_s\right) \right).
\end{split}
\end{equation*}
\qedsymbol

\section{Tools from the Literature}
We report here tools from the literature.

\begin{theorem}[{\citealp[Theorem~3.3.4]{yurinsky1995sums}}]
\label{theorem:bernstein-inequality-hilbert-space}
Let $Z_1, \dots, Z_n$ be independent random variables in a separable Hilbert space $\calH$ that are centered, that is,
\begin{equation*}
    \forall t \in [n], \E{Z_t} = 0,
\end{equation*}
If there exists real numbers $G, H \geq 0$ such that for any $p \geq 2$,
\begin{equation*}
    \sum_{t = 1}^{n} \bbE \nrm{Z_t}_{\calH}^p \leq \frac{1}{2} p! G^2 H^{p - 2},
\end{equation*}
it holds that for any $\eps > 0$,
\begin{equation*}
    \PR{\max_{1 \leq s \leq n }\nrm{   \sum_{t = 1}^{s} Z_t }_{\calH} > \eps} \leq 2 \exp \left( - \frac{\eps^2/2}{G^2 + \eps H} \right).
\end{equation*}
\end{theorem}

\begin{theorem}[{\citealp[Theorem~3]{sampson2000concentration}}]
\label{theorem:talagrand-type-concentration-phi-mixing-processes}
Let $X_1, \dots, X_n$ be a stationary $\phi$-mixing sequence.
For every $\eps > 0$,
\begin{equation*}
\begin{split}
    \PR{Z \geq \bbE{Z} + \eps} &\leq \exp \left( - \frac{1}{8 \nrm{\Gamma}_2} \min \set{ \frac{\eps}{C}, \frac{\eps^2}{4 \E{V}} } \right) \\
\end{split}
\end{equation*}
where 
\begin{equation*}
    Z \eqdef \sup_{g \in \calF} \abs{ \sum_{t=1}^{n} g(X_t) },
\end{equation*}
with $\calF$ an arbitrary class of real bounded functions with $\abs{g} \leq C$,
and where the random variable $V$ is defined as
\begin{equation*}
    V \eqdef \sum_{t=1}^{n} \sup_{g \in \calF} g(X_t)^2,
\end{equation*}
and $\Gamma$ is the coupling matrix (refer to Equation~\ref{definition:coupling-matrix}).
\end{theorem}

\end{document}